%% file: ms.tex
\documentclass[12pt]{article}
\usepackage{amsfonts,amssymb,amsmath}
\usepackage{a4,graphics,epsfig,psfrag}
\usepackage[utf8]{inputenc}
\usepackage{psfrag}
\usepackage{color}
\usepackage{tabularx}
\newcolumntype{Y}{>{\raggedright\arraybackslash}X}

\input{notations}
\begin{document}
\begin{center}
{\LARGE \bf 
Scaled Boundary Parametrizations \\[2mm]
in Isogeometric Analysis} \\[4mm]
Clarissa Arioli$^1$, Alexander Shamanskiy$^1$, Sven Klinkel$^2$, Bernd Simeon$^1$ \\[1mm]
$^1$ TU Kaiserslautern, Dept. of Mathematics \\
$^2$ RWTH Aachen, Chair of Structural Analysis
\end{center}

\begin{quote} \small
{\bf Abstract:}
This paper deals with a special class of parametrizations for Isogeometric Analysis (IGA). The so-called scaled boundary parametrizations are easy to 
construct and particularly attractive if only a boundary description of the computational domain is available. The idea goes back to the Scaled Boundary Finite Element Method (SB-FEM), which has recently been extended to IGA. We take here a different viewpoint and study these
parametrizations as bivariate or trivariate B-spline functions that are 
directly suitable for standard Galerkin-based IGA.
Our main results are first a general framework for this class of parametrizations, 
including aspects such as smoothness and regularity as well as generalizations to domains that are not star-shaped. Second, 
using the Poisson equation as example, we explain the relation 
between standard Galerkin-based IGA and the Scaled Boundary IGA by means of the Laplace-Beltrami operator. Further results concern 
the separation of integrals in both approaches and an analysis of the 
singularity in the scaling center. Among the computational examples 
we present a planar rotor geometry that stems from a screw compressor machine and compare different parametrization strategies.\\

\end{quote}

\input{introduction}

\input{parametrization}

\input{iga}

\input{singularity}

\input{computational} 

\section{Conclusions}

As we have seen, SB-parametrizations are easy to construct and 
can directly be cast into the framework of bi- or trivariate B-splines.
The methodology can be extended to include curved rays from the scaling center to the boundary, which relaxes the requirement of a star-shaped domain. Moreover, NURBS for the boundary description are suitable as well,
and the connection with art gallery problems indicates how more general geometries can be decomposed into star-shaped blocks.
SB-parametrizations are directly applicable in standard Galerkin-based IGA, 
but it turns out that the multiplicative structure of the Jacobian implies
a separation of the integrals, which can be exploited in order to speed up 
the computation of the stiffness matrix. The SB-IGA methods, on the other
hand, directly exploit this feature. In both cases, the singularity in the scaling center seems to be quite harmless, but this issue requires further investigation with respect to the overall convergence and to the condition number of the stiffness matrix.

\subsection*{Acknowledgement}
{\small We are grateful to Daniela Fu\ss eder, Carlo Lovadina and Angelos Mantzaflaris for fruitful discussions on the topic of 
scaled boundary parametrizations and for hints on related fields. 
Andreas Brümmer and Matthias Utri supplied us with the point cloud data for the rotor example in Fig.~\ref{fig:rotors}.
Moreover, we acknowledge the support by the European Union under grant no. 678727 (Alexander Shamanskiy, Project MOTOR) and by the German Research 
Council (DFG) under grant no. SI 756/5-1 (Clarissa Arioli, Project YASON). }

\input{appendix}

\bibliographystyle{plain}
\bibliography{igareferences}

\end{document}

%% file: notations.tex
\newcommand{\RR}{\mathbb R}

\newcommand{\f}[1]{\mathbf{#1}}

\newcommand{\fg}[1]{{\mbox{\boldmath $ #1 $}}}

%% file: introduction.tex
\section{Introduction}
Parametrizations of the computational domain play a crucial role in Isogeometric Analysis (IGA). We investigate here a special class of such parametrizations that can be viewed as generalizations of classical polar coordinates. 
These parametrizations can easily be constructed if only a boundary description is available. 
Our approach is inspired by the Scaled Boundary Finite Element Method (SB-FEM) \cite{Chen2015,Song1997,Song2000}
and its recent extension to IGA \cite{CHEN2016777,Klinkel2015,Natarajan2015}. 
We focus on standard Galerkin-based IGA in combination with such 
parametrizations and moreover 
address the connection to the Scaled Boundary IGA (SB-IGA).

Poisson's equation serves as model problem for our analysis. 
Making use of the framework of differential geometry, we express 
this partial differential equation by means of  the Laplace-Beltrami operator  and show that the weak form in parametric coordinates, which is at the core of IGA, can be derived in two ways. Either one derives first the weak form in physical coordinates and applies then the push-forward operator, or one applies first the coordinate transformation and then proceeds with the weak form. The first approach is taken in standard Galerkin-based IGA while the second is typical for the SB-IGA. The projection to a finite-dimensional spline space is then applied in both approaches, but numerical issues like the choice of quadrature rules will eventually lead to different implementations. In this context, we
derive an important structural property of scaled boundary parametrizations, which is the separation of the integrals that are computed for the entries in the stiffness matrix. This feature is a consequence of the multiplicative structure of the Jacobian of the geometry function. It comes up naturally in the SB-IGA but is hidden in the Galerkin-based IGA. Moreover, the separation property is a special case of  the recently introduced low rank tensor approximation technique by Mantzaflaris et al.~\cite{Mantzaflaris2015,Mantzaflaris2017}.

For scaled boundary parametrizations (SB-parametrizations) there are two drawbacks, which are the restriction to star-shaped domains and the singularity in the origin of the coordinate system. We address both issues and show how the methodology can be generalized. Most of our discussion concentrates on the planar case and standard tensor product B-splines, but we also comment on extensions 
to spatial geometries and NURBS (Non-Uniform Rational B-Splines).

In order to give a short overview on the state-of-the-art in the field, we refer
to the original work by Hughes et al.~\cite{Hughes2005,Cottrell.2009}
that set the ball rolling in IGA 
and to
Xu et al.~\cite{Xu2011} and Gravesen et al.~\cite{gravesen2014} for results on
the design and analysis of parametrizations in IGA. 
The SB-FEM has recently been extended towards crack propagation and other
nonlinear phenomena \cite{SAPUTRA2015213}. IGA for
partial differential equations on surfaces has been introduced by Dedè \& Quarteroni \cite{DEDE2015807}, exploiting also the powerful framework of differential geometry \cite{Berger2007}.
Singularities in IGA are addressed in Takacs \& J{ü}ttler \cite{TakacsJuettler2011,TakacsJuettler2012},
while very recently so-called multi-degree polar splines have been
proposed as a general alternative to circumvent the singularity in the scaling center by Toshniwal et al.~\cite{Toshniwal20171005}. 

The paper is organized as follows.
In Section 2, we define the class of SB-parametrizations and analyze their properties. Section 3 concentrates on their utilization in IGA, with emphasis on the Galerkin-based standard approach and side remarks on the relation to SB-IGA. Section 4 discusses the singularity in the scaling center both from an analytical and from a numerical perspective while Section 5 presents computational examples, among them a rotor geometry that stems from a screw compressor machine ~\cite{utri2017energy}.

%% file: parametrization.tex
\section{Scaled Boundary Parametrizations}

\begin{figure}
\centering
\includegraphics[height = 4cm]{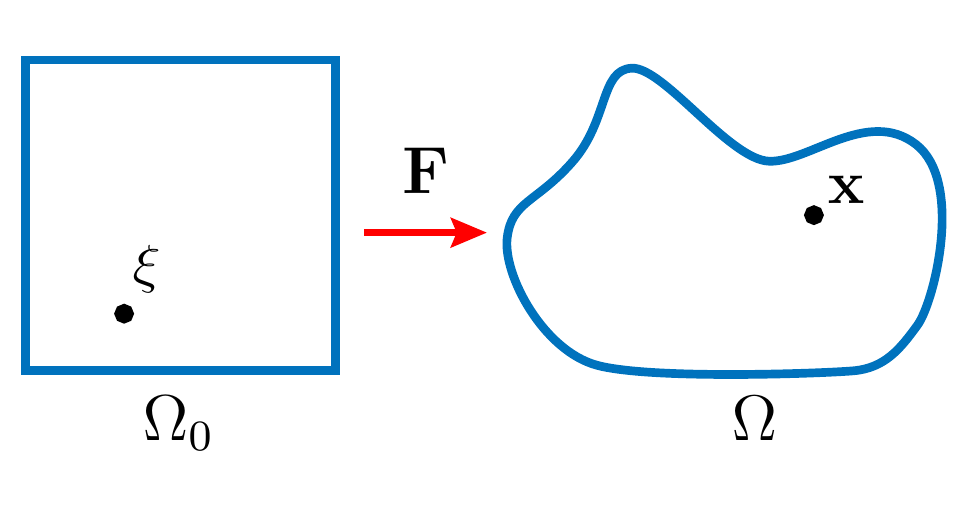}
\caption{Parametrization of $\Omega$. The parametric domain 
is denoted by $\Omega_0$.}\label{fig:geofun}
\end{figure}

Consider a domain $\Omega \subset \RR^d$ where $d=2$ or $d=3$ and its parameterization by a global geometry function 
\begin{equation}\label{num:geofun}
   \f{F}: \, \Omega_0 \rightarrow \Omega, \quad
             \f{F}(\fg{\xi}) = \f{x} = \left( \begin{array}{c} x_1 \\ \vdots \\ x_d
	                              \end{array}\right),
\end{equation}
see Fig.~\ref{fig:geofun}.
Below we will apply B-splines (or NURBS - Non Uniform Rational B-Splines) to define $\f{F}$, but for the moment the geometry function is simply an invertible $C^1$-mapping from the parameter domain $\Omega_0 \subset \RR^d$ to the physical domain. In our framework, $\Omega_0 = [0,1]^d$ is the unit square or unit cube, 
respectively.

Integrals over $\Omega$ can be transformed into integrals over $\Omega_0$ by means of the well-known integration rule
\begin{equation}\label{num:trafo}
       \int_\Omega w(\f{x}) \, \mathrm d\f x = 
       \int_{\Omega_0} w(\f{F}(\fg{\xi})) \,| \mbox{det } \f{DF}(\fg{\xi}) | \,
       \mathrm d\fg{\xi}
\end{equation}
with $d \times d$ Jacobian matrix $\f{DF}(\fg{\xi}) = \left( \partial F_i / \partial \xi_j \right)_{i,j=1:d}$.

In Computed-Aided Geometric Design (CAGD), objects are typically modelled in terms of their inner and outer hull, i.e., only a surface description is generated. In order to be able to apply IGA, however, one needs a computational mesh also in the interior that serves as discretization. 
The standard tool in IGA for this purpose are tensor product B-splines.

\subsection{Tensor Product B-splines}
To fix the notation, we shortly outline the construction of tensor product B-splines
in the case of $d=2$ dimensions. First we specify the 
polynomial degrees $(p,q)$ and the horizontal and vertical knot vectors 
\[
\Xi = \left\{\xi_1\le \xi_2\le\ldots \le \xi_{m+p+1}\right\},\quad 
\Psi = \left\{\eta_1\le \eta_2\le\ldots \le \eta_{n+q+1}\right\},
\]
which contain non-decreasing parametric real values so that
\[
0\le\mu(\Xi,\xi)\le p+1 \qquad \text{and} \quad
0\le\mu(\Psi,\eta)\le q+1 
\]
are the multiplicities of the parameter values in the knot vectors 
(the multiplicity $\mu(X,x)$ is zero if the given value $x$ is not a knot in $X$).
We write $M_{i,p}(\xi)$  and $N_{j,q}(\eta)$ for the univariate B-splines, 
which are computed by means of the recursion \cite{Boor1978}
\begin{align} 
N_{j,0}(\eta) &= \left\{ \begin{array}{ll}
1 & \text{for } \eta_j\le \eta < \eta_{j+1},\\
0 & \text{otherwise},
\end{array} \right.
\label{eq:bspline1}\\
N_{j,q}(\eta) &= 
\frac{\eta-\eta_j}{\eta_{j+q}-\eta_j}N_{j,q-1}(\eta) +
\frac{\eta_{j+q+1}-\eta}{\eta_{j+q+1}-\eta_{j+1}}N_{j+1,q-1}(\eta)
\label{eq:bspline2}
\end{align}
for $j=1,\ldots, n$. Fractions with zero denominators are considered zero. The same recursion applies to $M_{i,p}(\xi)$, replacing $j,q$ in 
(\ref{eq:bspline1})--(\ref{eq:bspline2}) by $i,p$ for $i=1,\ldots, m$. 

A function $\f{f}(\xi,\eta): 
[\xi_{p+1},\xi_{m+1}]\times[\eta_{q+1},\eta_{n+1}] \rightarrow \mathbb{R}^2$ is called a bivariate tensor product B-spline function if 
it has the form
\begin{equation}\label{eq:bsplinef}
\f{f}(\xi,\eta)=\sum_{i=1}^{m}\sum_{j=1}^{n}  M_{i,p}(\xi)N_{j,q}(\eta)  \f{d}_{i,j} 
\end{equation}
with the de Boor control points $\f{d}_{i,j} \in \mathbb{R}^2$, which form the control net
associated to the parametric representation. In order to work with the unit square as 
parameter domain, combined with open knot vectors, we assume
additionally $\xi_1=\ldots =\xi_{p+1} = 0$ and 
$\xi_{m+1} =\ldots = \xi_{m+p+1} = 1$ 
at the beginning and at the end; analogously  $\eta_1=\ldots =\eta_{q+1} = 0$ and 
$\eta_{n+1} = \eta_{n+2}=\ldots = \eta_{n+q+1} = 1$.

In $d=3$ dimensions, a trivariate B-spline tensor product function is generated
analogously, which results in 
\begin{equation}\label{eq:bsplinef3d}
\f{f}(\xi,\eta,\psi)=\sum_{i=1}^{m}\sum_{j=1}^{n}\sum_{k=1}^\ell  M_{i,p}(\xi)N_{j,q}(\eta)L_{k,r}(\psi)  \f{d}_{i,j,k} 
\end{equation}
with an additional set of univariate B-splines $L_{k,r}$ and control 
points $\f{d}_{i,j,k} \in \mathbb{R}^3$.

For the straightforward extension to NURBS, we refer to \cite{Piegl1997}.

\begin{figure}[b]
\centering
\includegraphics[height = 2.9cm]{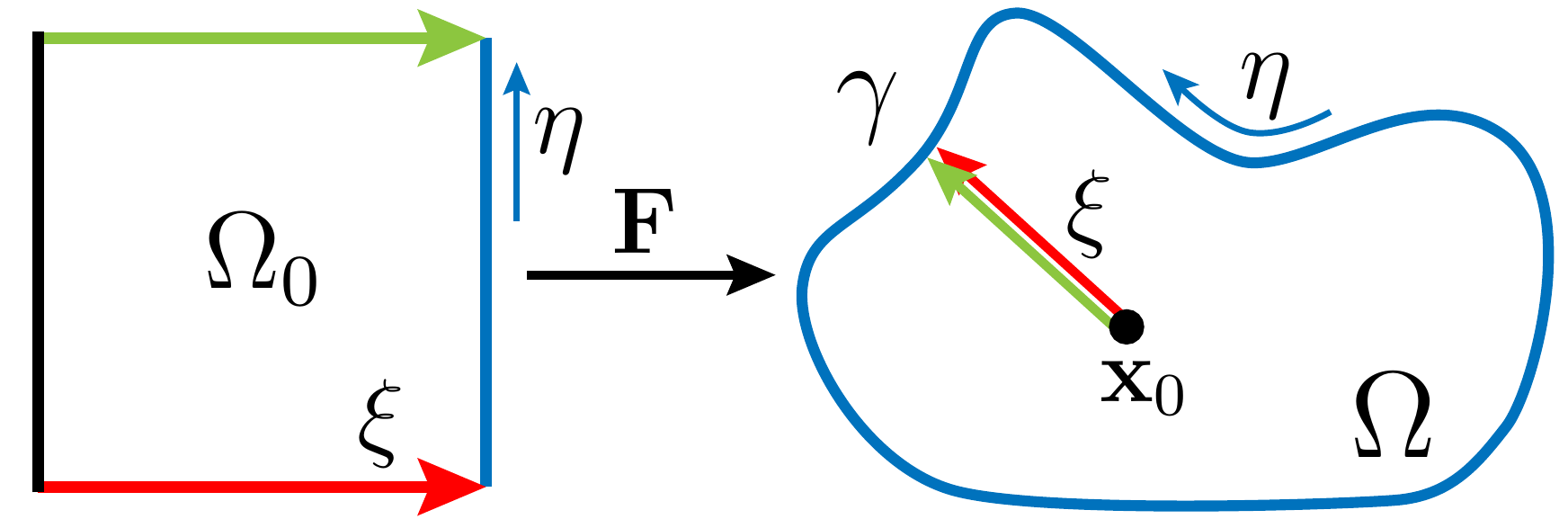}
\caption{Basic idea of an SB-parametrization.}\label{fig:basicidea}
\end{figure}

\subsection{The Basic Construction}
In $d=2$ dimensions, an SB-parametrization can be constructed as follows. We assume a star-shaped domain whose boundary is described by 
a spline curve $\fg{\gamma}$,
\begin{equation}\label{def:gamma}
         \fg{\gamma}(\eta) := \sum_{j=1}^{n} N_{j,q}(\eta) \f{c}_j
\end{equation}
with univariate B-splines $N_{j,q}$ of a certain degree $q$ and control points $\f{c}_j \in \mathbb{R}^2$. We require $\eta \in [0,1]$ for the independent variable and consider a closed curve with
$\fg{\gamma}(0) = \fg{\gamma}(1)$.  
This can be achieved by an open knot vector, as described in the previous subsection,  and control points 
$\f{c}_1 = \f{c}_n$. Other situations 
with several curve segments or wedge-shape geometries will be discussed below.

Next, we pick a point $\f{x}_0 \in \Omega$, the {\em scaling center}, and 
connect it with the boundary by means of rays that emanate from it, 
Fig.~\ref{fig:basicidea}. A ray that reaches from the scaling center to a control point $\f{c}_j$ on the boundary can be parametrized as
\begin{equation}\label{def:ray}
          (1-\xi)  \f{x}_0 + \xi \f{c}_j, \quad \xi \in [0,1]. 
\end{equation}
Replacing the control point $\f{c}_j$ by any point on the curve (\ref{def:gamma}), 
we obtain by the partition of unity
\begin{eqnarray}\label{def:geoF} 
          (1-\xi)  \f{x}_0 + \xi \fg{\gamma}(\eta) &=& 
         (1-\xi)  \sum_{j=1}^n N_{j,q}(\eta) \f{x}_0 + \xi \sum_{j=1}^{n} N_{j,q}(\eta) \f{c}_j \nonumber      \\
                  & = & \sum_{i=1}^{2}\sum_{j=1}^n M_{i,1}(\xi) N_{j,q}(\eta) \f{d}_{ij}     =:         \f{F}(\xi,\eta)  \label{def:geoF}  
\end{eqnarray}
 with control points $\f{d}_{1j} := \f{x}_0,\,                                          \f{d}_{2j} := \f{c}_j, \, j=1,\ldots,n;$
and linear B-splines $M_{1,1}(\xi) = 1-\xi, \, \, M_{2,1}(\xi) = \xi$.

The bivariate parametrization $\f{F}(\xi,\eta)$ maps the unit square 
$\Omega_0 = [0,1]^2$ to the domain $\Omega$ with a multiple control point
$\f{d}_{1j} = \f{x}_0$ in the scaling center and the other control points 
$\f{d}_{2j} = \f{c}_j$ adopted from the boundary curve. In other words,
the left vertical edge of the unit square collapses to the scaling center
while the right vertical edge is mapped to the boundary, 
Fig.~\ref{fig:basicidea}.
Note that the geometry function $\f{F}$ in (\ref{def:geoF}) features the usual tensor product structure (\ref{eq:bsplinef})
with only $m=2$ linear B-splines in $\xi$-direction and $n$ B-splines inherited from the boundary curve in $\eta$-direction.

In some cases it will be advantageous to replace the double sum in (\ref{def:geoF}) 
by the compact form
\begin{equation}\label{def:geoFc}
    \f{F}(\xi,\eta)  = \f{x}_0 + \xi (\f{C} \cdot \f{N}(\eta)  - \f{x}_0),
\end{equation}
which is the usual notation in the SB-FEM 
\cite{Song1997,Song2000} and SB-IGA\cite{Chen2015,CHEN2016777,Natarajan2015}. Here, the matrix
$\f{C} := (\f{c}_1, \ldots, \f{c}_n) \in \mathbb{R}^{2\times n}$ 
contains the control points on the boundary and the vector
$\f{N}(\eta) := (N_{1,q}(\eta), \ldots, N_{n,q}(\eta))^T \in \mathbb{R}^n$
the B-splines.

Finally, we extend the bivariate spline parametrization (\ref{def:geoF}) 
and apply both knot insertion and degree elevation in radial direction,
 which preserves the rays that emanate from the scaling center but
 leads to a more general formulation
 \begin{equation}\label{def:geoFc2}
\f{F}(\xi,\eta) =
 \sum_{i=1}^{m}\sum_{j=1}^n M_{i,p}(\xi) N_{j,q}(\eta) \tilde{\f{d}}_{ij} .
 \end{equation}
Note that the multiple control point in the scaling center is still present
here, i.e., $\tilde{\f{d}}_{1j} = \f{x}_0$ for $j=1,\ldots,n$ and that 
$\tilde{\f{d}}_{mj} = \f{c}_j$ for
the control points of the boundary curve, $j=1,\ldots,n$. 
The other extra control points $\tilde{\f{d}}_{ij} $ result from the 
usual steps of knot insertion and degree elevation. 
In this way, we obtain a finer parametrization that still preserves the 
scaled boundary idea. At the same time, it possesses a discretization 
in both $\xi$ and $\eta$ that can be used as natural isogeometric mesh for 
a numerical simulation.

\subsection{Properties and Extensions}

We aim now at studying the regularity and smoothness of the 
SB-parametrization. For this purpose, it is convenient to
use the representation (\ref{def:geoFc}), whose Jacobian reads
\begin{equation}\label{eq:DF}
\f{DF}(\xi,\eta) = \big( \f{C} \f{N}(\eta) - \f{x}_0 \, |  \, \xi \f{C} \f{N}^\prime(\eta)
\big) = 
\begin{pmatrix}
\f{c}_1^T\f{N}(\eta)-x_{0,1}  & \f{c}_1^T\f{N}^\prime(\eta) \\
\f{c}_2^T\f{N}(\eta)-x_{0,2}  & \f{c}_2^T\f{N}^\prime(\eta) 
\end{pmatrix}
\begin{pmatrix}
1 & 0\\
0&\xi
\end{pmatrix}.
\end{equation}
Here, we used the notation $\f{c}_1^T := \f{C}(1,:)$ for the first
and $\f{c}_2^T := \f{C}(2,:)$ for the second row of $\f{C}$, respectively.
Furthermore, $\f{x}_0^T = (x_{0,1},x_{0,2})$ and the derivative
with respect to $\eta$ is written as $\f{N}^\prime = \partial \f{N} / \partial \eta$.

The determinant of the Jacobian is
given by
\begin{equation}\label{eq:detDF}
 \textrm{det} \, \f{DF}(\xi,\eta) = \xi J(\eta), 
 \end{equation}
 where
 \[
 \quad J(\eta):=
\f{c}_1^T  \f{N}(\eta)\f{c}_2^T  \f{N}'(\eta) - x_{0,1}\f{c}_2^T  \f{N}'(\eta) 
- \f{c}_1^T  \f{N}'(\eta)\f{c}_2^T  \f{N}(\eta) + x_{0,2}\f{c}_1^T  \f{N}'(\eta).
\]
The multiplicative structure of the Jacobian allows a straightforward analysis of the
parametrization in terms of regularity and smoothness. 
Obviously, in the scaling center $\f{x}_0$ where $\xi = 0$ the parametrization $\f{F}$
becomes singular. For the derivative in $\eta$-direction, we have additionally
 $ \partial \f{F}(\xi=0,\eta) / \partial \eta = \f{0}$. 
 The smoothness of the boundary curve $\fg{\gamma}$, on the other hand,
 determines the smoothness of $\f{F}$: if $\gamma$ is of class $C^k$, then
 $\f{F}$ is of class $C^k$ as well. If $\gamma$ possesses points of reduced smoothness, e.g. an interpolatory knot where the control point $\f{c}_j$ lies
 on the curve, then the ray $(1-\xi)  \f{x}_0 + \xi \f{c}_j, \,\xi \in [0,1]$, which runs
 from the scaling center to $\f{c}_j$, will inherit the smoothness and form 
 a $C^0$-edge in the interior. This reasoning applies also to the regularity.
 If $\fg{\gamma}$ is injective, the parametrization will also be injective
 except for the scaling center.

In passing we note that recently so-called multi-degree polar splines have been
introduced as a general alternative to circumvent the singularity in the scaling center
\cite{Toshniwal20171005}. 

We now turn to extensions of SB-parametrizations.

\begin{figure}
	\centering
	\includegraphics[height = 3.3cm]{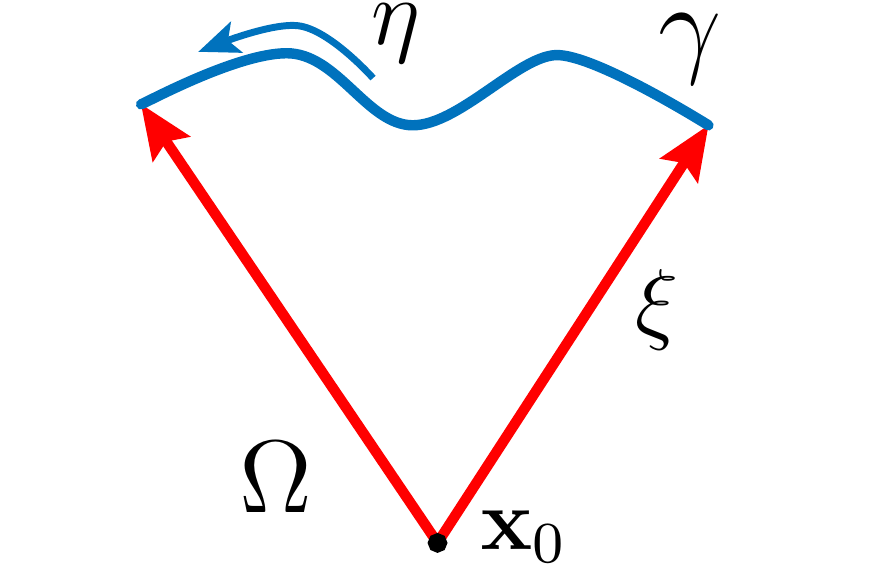} \mbox{\qquad}
	\includegraphics[height = 5cm]{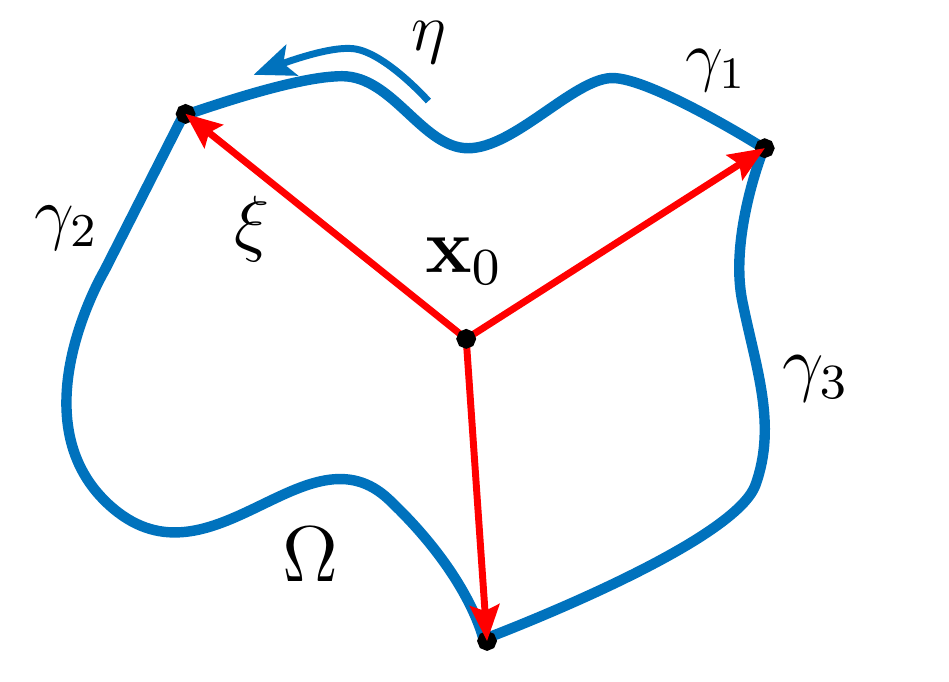}
	\caption{Wedge-shape domain $\Omega$ and combination of several wedges }\label{fig:wedge}
\end{figure}

{\bf Several Curve Segments.\/} There are situations where several curve segments 
are connected in order to define the computational domain $\Omega$. In this case,
one can decompose the domain into several wedge-shape objects where 
each wedge forms a triangle as sketched in 
Fig.~\ref{fig:wedge} on the left. The scaling center of such a wedge is in one of the vertices,
and the two edges that emanate from it are the rays that point to the end points
of the curve $\fg{\gamma}$. The representations (\ref{def:geoFc}) and
(\ref{def:geoFc2}) hold also for the parametrization of the wedge. If we
combine several curves as in Fig.~\ref{fig:wedge} on the right, 
we first pick a common scaling center $\f{x}_0$, and each 
wedge possesses then its own parametrization. 
Due to the connectivity of the 
boundary curve, however, it is guaranteed that the individual wedges do not intersect and that the domain $\Omega$ is completely covered.

One might wonder whether the situation of several curve segments is actually of relevance since it is always possible to define a new boundary curve that comprises the individual segments. E.g., in Fig.~\ref{fig:wedge} we can set
$\tilde{\fg{\gamma}}:= \fg{\gamma}_1 \cup \fg{\gamma}_2 \cup\fg{\gamma}_3$.
The new curve $\tilde{\fg{\gamma}}$ then might feature interpolatory knots with
reduced continuity, but this would mean no obstacle to use it as starting point for a single SB-parametrization (\ref{def:geoFc2}). 
Nevertheless, several curve segments are useful as they present a natural decomposition of the domain, similar to a multi-block or multi-patch strcuture, that
is inherited from the boundary geometry. Morever, this description provides insight
into the 3D case where each curve segment turns into a surface patch with an individual parametrization.

{\bf Extensions beyond Star-shaped Domains.\/} 
Obviously, the restriction to a star-shaped domain is a limitation of SB-parametrizations. We shortly discuss two ways to overcome this drawback.

\begin{figure}[h!]
	\centering
	\includegraphics[height = 7.13cm]{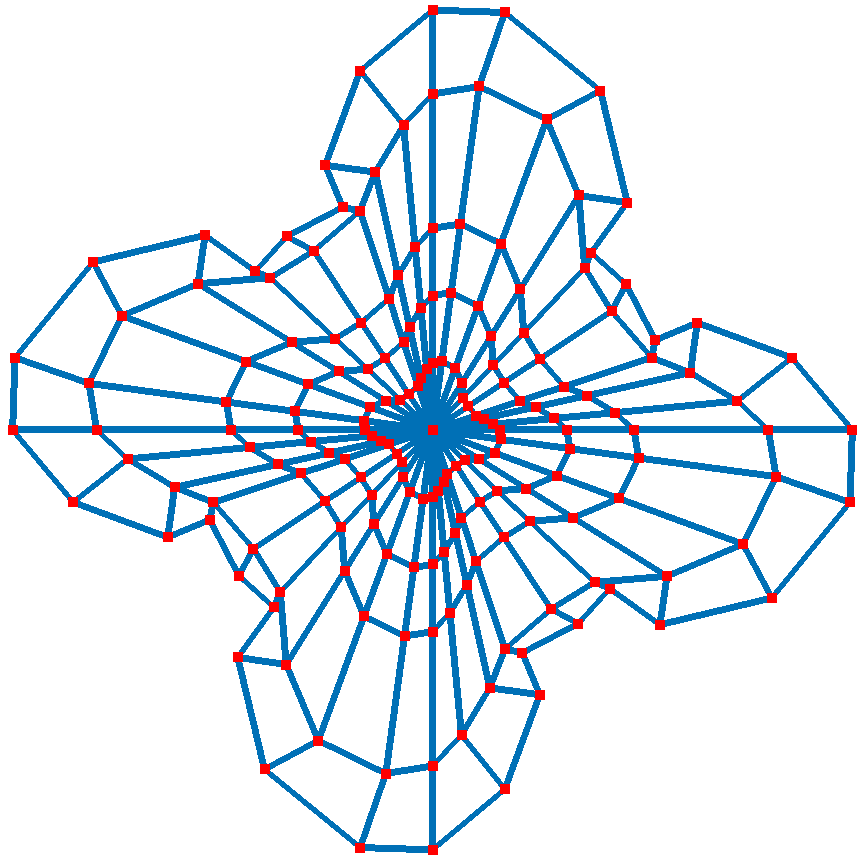}
	\includegraphics[height = 7.13cm]{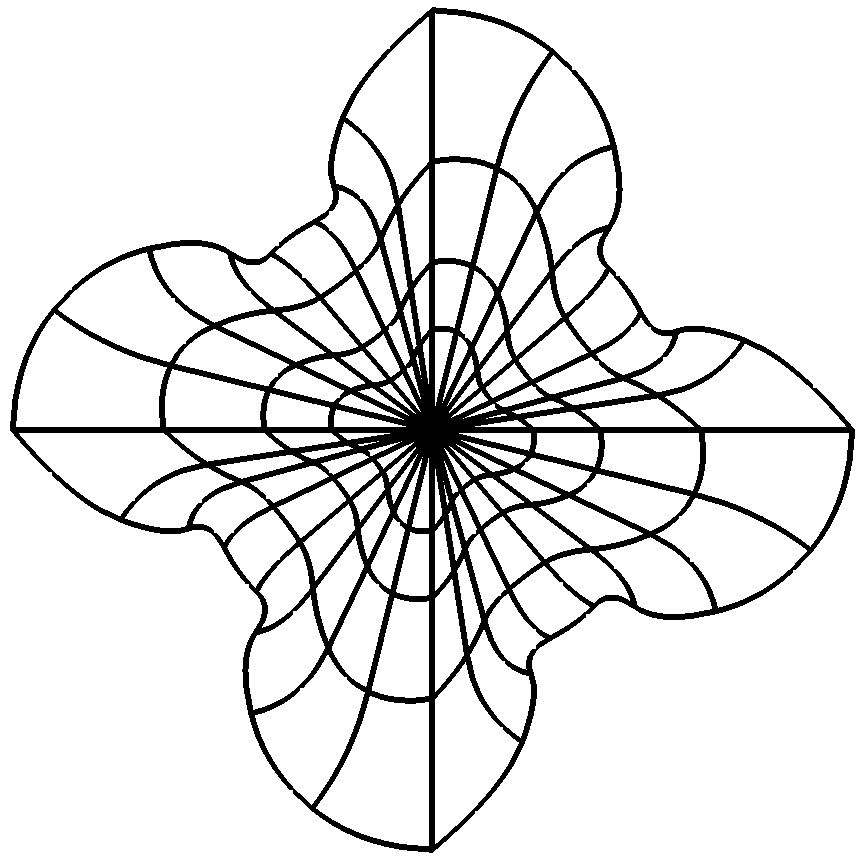}
	\caption{Screw compressor rotor as an example of a non-star-shaped domain. Control net with control points (left) and isoparametric lines (right).}
	\label{fig:rotors}
\end{figure}

Consider again the general parametrization (\ref{def:geoFc2}) where 
the control points in the interior are given by 
 $\tilde{\f{d}}_{ij}$ for $i=2,\ldots,m-1$ and $j=1,\ldots,n$.
If we fix an index $j$ and consider
the corresponding ray from the scaling center to the control point $\f{c}_j$ on the
boundary, it is possible to move the control points 
$\tilde{\f{d}}_{2,j}, \ldots, \tilde{\f{d}}_{m-1,j}$
along the ray such that we obtain a curved ray.
This technique can be used to relax the requirement of a star-shaped domain , see Fig.~\ref{fig:rotors} for an example with a rotor geometry that is part of a screw compressor unit. In this figure, the visible mesh on the right consists of the isoparametric lines $\xi = \xi_k$ and $\eta=\eta_\ell$, which must be distinguished from the control net on the left.

\begin{figure}
	\centering
	\includegraphics[height = 5cm]{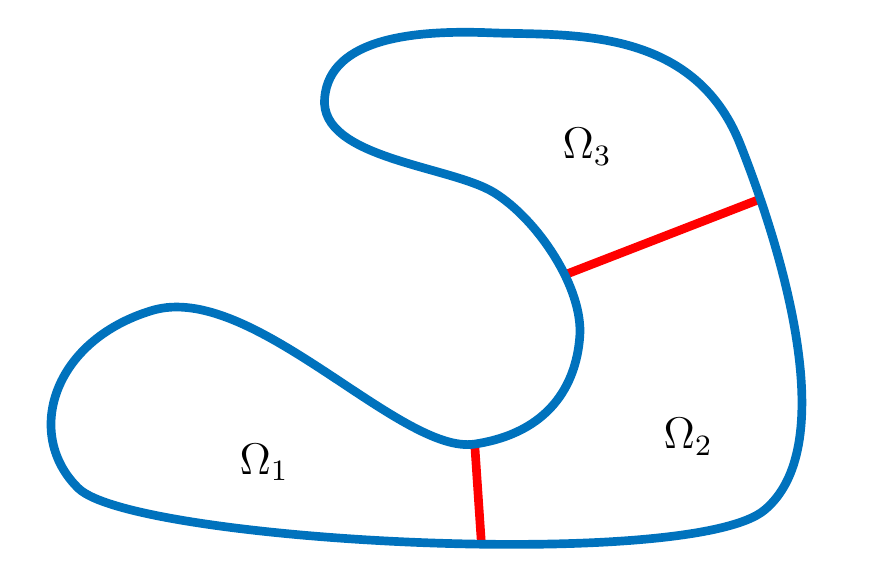}
	\caption{Non-star shape case of $\Omega$. }\label{fig:artgallery}
\end{figure}

A more general approach relies on methods for the so-called {\em art gallery problem}. Basically, this means that for a general domain $\Omega$ such as in 
Fig.~\ref{fig:artgallery}, we search for a decomposition into blocks or patches that 
are star-shaped. If the domain $\Omega$ has a polygonal boundary with $k$ vertices,
there is an upper bound on the number of blocks that are required for such 
a decomposition.  The bound is $\lfloor k/3 \rfloor$, and there are algorithms that 
compute the decomposition automatically \cite{ORourke1987}.

{\bf Extension to NURBS.\/} 
If the curve $\fg{\gamma}$ in (\ref{def:gamma}) is given as linear
combination of NURBS \cite{Piegl1997}
\begin{equation}\label{def:NURBScurve}
\fg{\gamma}(\eta) =  \sum_{j=1}^{n} R_{j,q}(\eta) \f{c}_j, \quad
                   R_{j,q}(\eta) = \frac{N_{j,q}(\eta) w_j}{\sum_{i=1}^n N_{i,q}(\eta) w_i}
\end{equation}
with weights $w_j, \, j=1,\ldots,n$, the above construction of the SB-parametrization stays the same. Moreover, all the properties discussed so far, including the multiplicative structure of the Jacobian, are preserved.

\begin{figure}[h]
	\centering
	\includegraphics[height = 7.13cm]{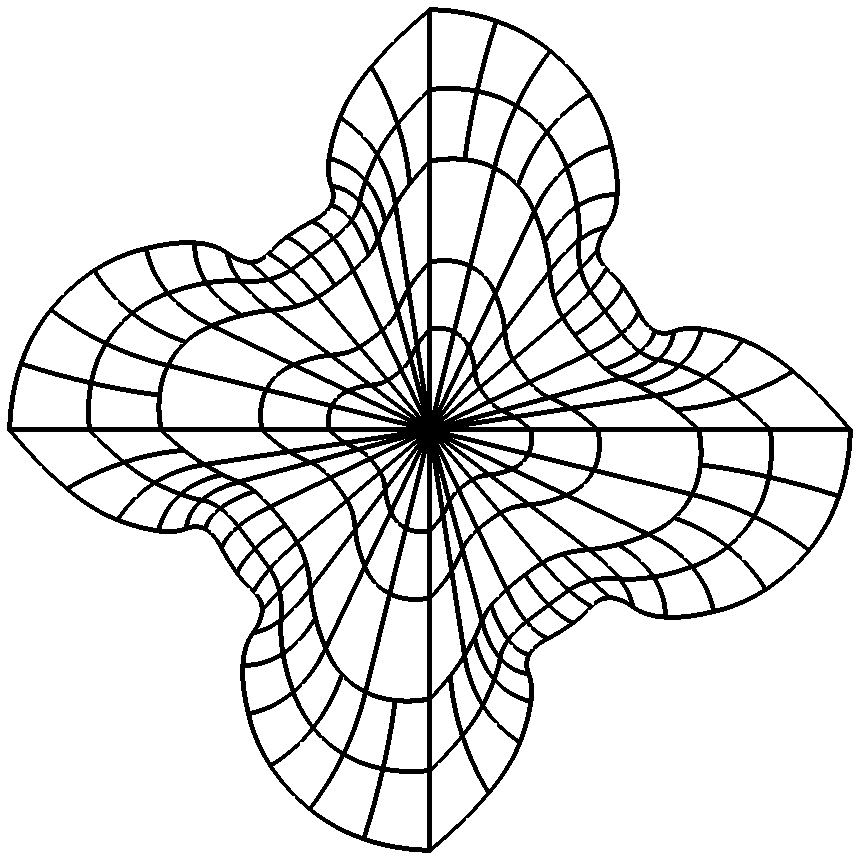}
	\caption{Screw compressor rotor. Local refinement along the boundary using THB-splines.}\label{fig:rotorTHB}
\end{figure}

{\bf Extension to THB.\/} 
It is possible to use local refinement along the boundary, where higher numerical accuracy is often necessary. Fig.~\ref{fig:rotorTHB} shows an example of local refinement using THB-splines \cite{giannelli2012thb,vuong2011hierarchical}.

{\bf Extension to 3D Geometries.\/}
Obviously, the basic idea for constructing an SB-parametrization also carries over to 3D objects that are given in terms of a surface description. For a detailed discussion of this case we refer to
\cite{CHEN2016777} where such parametrizations are constructed and
combined with SB-IGA. We shortly consider a star-shaped solid 
$\Omega$ with 
surface $\partial \Omega$ that is parametrized via 
\begin{equation}\label{def:sigma}
         \fg{\sigma}(\eta,\psi) := \sum_{j=1}^{n}\sum_{k=1}^\ell N_{j,q}(\eta)L_{k,r}(\psi) \f{c}_{jk}
\end{equation}
with bivariate B-splines $N_{j,q}L_{k,r}$ and control points $\f{c}_{jk} \in \mathbb{R}^3$. In general, there will be several surface patches 
that are glued together to form the complete surface, but we skip this 
technical detail. Picking a scaling center $\f{x}_0 \in \Omega$, 
we can proceed as in the planar case and 
write a ray from the center to the surface as
\begin{equation}\label{def:sb3d}
(1-\xi)  \f{x}_0 + \xi \fg{\sigma}(\eta,\psi)
                  =  \sum_{i=1}^{2}\sum_{j=1}^n \sum_{k=1}^\ell
                  M_{i,1}(\xi) N_{j,q}(\eta) L_{k,r}(\psi) \f{d}_{ijk}     =:         \f{F}(\xi,\eta,\psi)   
\end{equation}
 with control points $\f{d}_{1jk} := \f{x}_0,\,                                          \f{d}_{2jk} := \f{c}_{jk}, \, j=1,\ldots,n; \, k=1,\ldots,\ell$.

For some geometries, one can also apply the scaled boundary idea
only in a cross-section and then extrude the geometry in the third direction, generalizing  the concept of cylindrical coordinates.

%% file: iga.tex
\section{Utilization in Isogeometric Analysis}

In this section we analyze the usage of SB-parametrizations in
IGA. 
For ease of presentation, we consider Poisson's equation
\begin{equation}\label{num:poisson}
     -\Delta u = f  \quad \mbox{in } \Omega
\end{equation}
as an illustrative model problem.
Here,  $\Omega \subset \RR^d$ is a domain with boundary 
$\partial \Omega$,
the function $f: \Omega \rightarrow \RR$ is a given source term, and the unknown function $u: \Omega \rightarrow \RR$ shall satisfy the Dirichlet boundary condition
\begin{equation}\label{num:dirichlet}
     u = 0 \quad \mbox{on } \partial \Omega .
\end{equation}
The discussion of non-zero Dirichlet and Neumann boundary conditions is
omitted to keep the exposition concise.
The weak form of the PDE (\ref{num:poisson}) is obtained by multiplication with test functions $v$ and integration over $\Omega$.
More specifically, one defines the function space
\begin{equation}\label{num:defV}
   V := \{ v \in H^1(\Omega), \, v = 0 \,\, \mbox{on } \partial \Omega \},
\end{equation}  
which consists of all functions $v \in L_2(\Omega)$ that possess weak and square-integrable first derivatives and that vanish on the boundary. 
For functions $u,v \in V$, the bilinear form
\begin{equation}\label{num:bilinear}
    a(u,v) := \int_\Omega \nabla u \cdot \nabla v \, \mathrm d\f x
\end{equation}
is well-defined, and even more, it is symmetric and coercive. 
Setting  
\begin{equation}\label{num:linform}
    \langle l, v \rangle := \int_\Omega f v \, \mathrm d \f{x}
\end{equation}
as linear form for the integration of the right hand side, the solution $u \in V \subset H^1(\Omega)$ is then characterized by the weak form
\begin{equation}\label{num:weakform}
     a(u,v) = \langle l, v \rangle  \quad \mbox{for all } v \in V 
\end{equation}
and the boundary condition $u = 0$ (in the sense of traces).

\subsection{Transformation of the Problem Formulation}

Next assume that we have a parametrization of $\Omega$ available as in
(\ref{num:geofun}).
For the differentiation, the chain rule applied to $u(\f{x}) = u(\f{F}(\fg{\xi}))
=:\hat{u}(\fg{\xi}) $ yields, using a row vector notation for the gradient $\nabla u$, 
\begin{equation}\label{num:chainrule}
     \nabla_{\f{x}} \, u(\f{x}) = \nabla_{\f{\xi}}\, \hat{u}(\fg{\xi}) 
                           \cdot \f{DF}(\fg{\xi})^{-1} .
\end{equation}
The integrals in the weak form (\ref{num:weakform}) then satisfy the transformation rules
\begin{equation}\label{num:trafoa}
    \int_\Omega \nabla u \cdot \nabla v \, \mathrm d\f x  
    = \int_{\Omega_0}  (\nabla \hat{u} \, \f{DF}^{-1}) \cdot 
    (\nabla \hat{v} \, \f{DF}^{-1}) \,
    |\mbox{det } \f{DF} | \,
       \mathrm d\fg{\xi} 
\end{equation}
and 
\begin{equation}\label{num:trafob}
\int_\Omega f v \, \mathrm d\f x = 
       \int_{\Omega_0} \hat{f} \hat{v} \,| \mbox{det } \f{DF} | \,
       \mathrm d\fg{\xi} \,.
\end{equation}
The integrals over $\Omega_0$ that are defined by the right hand sides of
 (\ref{num:trafoa}) and (\ref{num:trafob}) are at the core of Galerkin-based IGA. 
 
 For our purposes, it is now advantageous to adopt a more general 
 viewpoint from differential geometry where the Laplace operator 
 is extended to operate on functions defined on surfaces in Euclidean space and on Riemannian manifolds, see Berger \cite{Berger2007}. We define the metric tensor or first fundamental form
 \begin{equation}\label{eq:defgmetric}
  \f{g}(\fg{\xi}) := \f{DF}(\fg{\xi})^T\f{DF}(\fg{\xi}),
 \end{equation}
 which has the determinant $\mbox{det } \f{g} = (\mbox{det } \f{DF})^2$.
 The transformed weak form from  
 (\ref{num:trafoa}) and  (\ref{num:trafob}) 
 then reads
 \begin{equation}\label{num:trafometric}
    \int_{\Omega_0}  \nabla \hat{u} \, \f{g}^{-1} \cdot
    \nabla \hat{v} 
    |\mbox{det } \f{g} |^{1/2} \,
       \mathrm d\fg{\xi}  = \int_{\Omega_0} \hat{f} \hat{v} \,| \mbox{det } \f{g} |^{1/2} \,
       \mathrm d\fg{\xi}
\end{equation}
for test functions $\hat{v} = v \circ \f{F}$. In terms of regularity and smoothness, the parametrization $\f{F}$ is here required to be of class
$C^1$ and to be injective almost everywhere, which means that singularities
in sets of measure zero, such as points, are admitted.

Corresponding to the weak form (\ref{num:trafometric}), there 
is a strong form that makes use of the Laplace--Beltrami operator. 
Like the standard Laplacian, this operator is defined as the divergence of the gradient in parametric coordinates, and it reads
\begin{equation}\label{eq:Beltrami}
  \Delta_\xi \hat{u} = |\mbox{det } \f{g} |^{-1/2}
  \sum_{k=1}^d \frac{\partial}{\partial \xi_k}
  \left( |\mbox{det } \f{g} |^{1/2} \nabla \hat{u} \, \f{g}^{-1} 
  \right).
\end{equation}
For the strong form of Poisson's equation with respect to the metric 
$\f{g}$ and expressed in the coordinates $\fg{\xi}$, this implies the representation
\begin{equation}\label{eq:strongxi}
-   \sum_{k=1}^d \frac{\partial}{\partial \xi_k}
  \left( |\mbox{det } \f{g} |^{1/2} \nabla \hat{u} \, \f{g}^{-1} 
  \right)
= |\mbox{det } \f{g} |^{1/2} \hat{f}   \quad \mbox{in } \Omega_0.
\end{equation}
Clearly, this requires now the parametrization $\f{F}$ to be of class $C^2$ and 
injective in order to hold in all of $\Omega_0$.
\subsubsection*{Equivalence of the Two Approaches}
Based on this general framework provided by differential geometry,
there are two ways to derive the weak form (\ref{num:trafometric}) in parametric coordinates $\fg{\xi}$ and to proceed then with a discretization. We start with the strong form (\ref{num:poisson}) in Cartesian coordinates $\f{x}$ and then have the
following choices:  
\begin{itemize}
\item[(i)] First we derive the weak form (\ref{num:weakform}) in $\f{x}$
and then we apply the transformation to the weak form  (\ref{num:trafometric}) in $\fg{\xi}$, followed by the projection
to a finite-dimensional spline space.
\item[(ii)] First we apply the transformation to the strong form 
       (\ref{eq:strongxi}) in $\fg{\xi}$ and then we proceed with the weak form 
       (\ref{num:trafometric}) and the projection to a finite-dimensional spline space. 
\end{itemize}
Option (i) is the approach taken in the Galerkin-based IGA while
option (ii) is the route taken in the Scaled Boundary IGA. In the latter, however, it is not mandatory to use the weak form  (\ref{num:trafometric}) 
in all independent variables. Instead,
the numerical treatment is typically split into a discretization of the boundary
and a subsequent procedure for the radial direction, see Subsection \ref{subs:scaledboundIGA} below.
From a theoretical point of
view, both (i) and (ii) lead to the same weak form in parametric
coordinates.

To conclude this digression on manifolds, we remark that in a sense IGA for the Poisson equation (\ref{num:poisson}) 
can be interpreted as solving the transformed PDE  (\ref{eq:strongxi})
with the Laplace-Beltrami operator on the unit square or unit cube.

As a classical example, we recall polar coordinates for the unit disk $\Omega\subset \mathbb{R}^2$, which read
 $$\f{F}(\xi,\eta)=\xi \begin{pmatrix}
\cos\eta\\
\sin\eta
\end{pmatrix}, \quad \xi \in [0,1], \,\, \eta \in [0, 2\pi],$$
and result in the metric tensor 
\[   \f{g}(\xi,\eta) = \left( \begin{array}{cc} 1 & 0 \\ 0 & \xi^2 \end{array} \right).
\]
Then the transformed weak form (\ref{num:trafometric}) becomes
\begin{equation}\label{num:weakandpolar}
\int_{\Omega_0}{\left(\xi\hat{u}_\xi\hat{v}_\xi+\frac{1}{\xi}\hat{u}_\eta\hat{v}_\eta \right)\mathrm d{\xi}\mathrm d{\eta}}
 = \int_{\Omega_0} \xi \hat{f}\hat{v} d{\xi}\mathrm d{\eta},
\end{equation}
where $\hat{u}_\xi=\partial \hat{u}/\partial \xi$ and so forth.

Alternatively to (\ref{num:weakandpolar}), the strong form 
of the Poisson problem is given by
\begin{equation}\label{eq:polarstrong}
 -\left(\hat{u}_{\xi\xi} + \frac{1}{\xi} \hat{u}_\xi
 +\frac{1}{\xi^2}\hat{u}_{\eta\eta}\right) = \hat{f} .
 \end{equation}
Note the Laplace operator in polar coordinates on the left hand side.
This example gives valuable insight into the structure of the transformed
PDE. For the more general case of an SB-parametrization, 
this structure is similar with factors $\xi^{-1}$ and $\xi^{-2}$ that 
indicate the singularity in the scaling center, as we will see next.

\subsection{Galerkin-based IGA}
The Galerkin projection of the weak form (\ref{num:weakform}) 
in physical coordinates $\f{x}$
replaces the infinite dimensional space $V$ by a finite dimensional subspace $V_h \subset V$, with the subscript $h$ indicating the relation to a spatial grid. 
Let $\phi_1, \ldots, \phi_K$ be a basis of $V_h$, then the numerical approximation $u_h$ is constructed as linear combination
\begin{equation}\label{num:defuh}
     u_h =  \sum_{i=1}^{K} q_i \phi_i
\end{equation}
with unknown real coefficients $q_i \in \RR$. 

Upon inserting $u_h$ into the weak form (\ref{num:weakform}) and testing with $v = \phi_j$ for $j=1,\ldots,K$, one obtains the linear system
\begin{equation}\label{num:linsys}
    \f{A} \f{q} = \f{r}
\end{equation}
with $K \times K$ stiffness matrix $\f{A} = ( a(\phi_i,\phi_j))_{i,j=1,\ldots,K}$ and right hand side vector $\f{r} = ( \langle l, \phi_i \rangle)_{i=1,\ldots,K}$.
Since the matrix $\f{A}$ inherits the properties of the bilinear form $a$, it is straightforward to show that $\f{A}$ is symmetric positive definite, and thus the numerical solution $\f{q}$ or $u_h$, respectively, is well-defined.

In Galerkin-based IGA, we define 
 $\phi_i = \psi_i \circ \f{F}^{-1}, \, i=1,\ldots,K,$  
as basis functions
via the push forward operator, 
 where $\psi_i = M_{j,p}N_{k,q}$ is 
given by bivariate B-splines in the planar case
and $\psi_i = M_{j,p}N_{k,q}L_{w,r} $ by trivariate B-splines in the spatial case.
The total number of degrees of freedom is accordingly $K = mn$ or
$K = mn\ell$, respectively, and the numerical approximation 
becomes, expressed in parametric coordinates, 
a bivariate B-spline tensor product function (\ref{eq:bsplinef})
or a trivariate function (\ref{eq:bsplinef3d}).
Overall, this projection step then boils down to 
inserting $ \hat{u} = \sum q_i \psi_i$ and $\hat{v} = \psi_j, j=1,\ldots,K,$ into the  weak form (\ref{num:trafometric}).

\subsubsection*{Exploiting the Structure of the Scaled Boundary Parametrization}

If the parametrization of the computational domain is given in terms 
of tensor product B-splines, any IGA solver can be applied directly. 
But such a black box procedure misses an important structural feature
of SB-parametrizations.

To illustrate this, we consider the domain $\Omega$ in the case $d=2$, parametrized by the geometry function $\f{F}$ as given in (\ref{def:geoF}) and (\ref{def:geoFc}). 

For the inverse of the metric tensor $\f{g}$ one obtains
\begin{align*}\renewcommand*{\arraystretch}{1.8}
 \f{DF} ^{-1}(\xi,\eta)  \f{DF} ^{-T}(\xi,\eta)=& \frac{1}{(\mbox{det } \f{DF})^2}
 \begin{pmatrix}
 \xi & 0 \\
 0 &1 
 \end{pmatrix}
\left(\begin{matrix}
\f{b}_1^T(\eta)\\
\f{b}_2^T(\eta)
\end{matrix}\right)
\left(\begin{matrix}
	\f{b}_1(\eta) &|\quad	\f{b}_2(\eta)
\end{matrix}\right)
\begin{pmatrix}
\xi & 0 \\
0 &1 
\end{pmatrix}\\[1mm]
=&\dfrac{1}{ J(\eta)^2}\left(\begin{matrix}
\f{b}_1^T(\eta)\f{b}_1(\eta) & \f{b}_1^T(\eta)\f{b}_2(\eta)/\xi\\
\f{b}_2^T(\eta)\f{b}_1(\eta)/\xi & \f{b}_2^T(\eta)\f{b}_2(\eta)/\xi^2
\end{matrix}\right)
\end{align*}
with \begin{equation*}
\f{b}_1^T(\eta) := (\f{c}_2^T\f{N}^\prime(\eta),\,-\f{c}_1^T\f{N}^\prime(\eta)),\quad
\f{b}_2^T(\eta) := (-\f{c}_2^T\f{N}(\eta)+x_{0,2},\,\f{c}_1^T\f{N}(\eta)-x_{0,1}).
\end{equation*}
Then (\ref{num:trafoa}) becomes, omitting the arguments $\xi,\eta$ for a compact notation,
\begin{eqnarray}
\int_{\Omega_0}  \nabla \hat{v}^T \, \f{DF}^{-1}
    \, \f{DF}^{-T} \nabla \hat{u} \,
    |\mbox{det } \f{DF} | \,
       \mathrm d\xi \mathrm d \eta \mbox{\hspace{6cm}} & &
\nonumber \\ 
	= \int_0^1\int_0^1 {\small \frac{1}{J}}\left(\xi \,\hat{u}_\xi\hat{v}_{\xi} \,\f{b}_1^T\f{b}_1+\hat{u}_\eta\hat{v}_\xi \,\f{b}_1^T\f{b}_2+\hat{u}_\xi\hat{v}_\eta \,\f{b}_2^T\f{b}_1
	+ \frac{1}{\xi}\hat{u}_\eta\hat{v}_\eta \,\f{b}_2^T\f{b}_2
	\right) \;\text{d}\xi\text{d}\eta .  
	\label{num:weakwithspecF}
\end{eqnarray}
The transformed bilinear form (\ref{num:weakwithspecF}) 
exhibits the underlying structure of the SB-parametrization.
For the numerical solution, (\ref{num:weakwithspecF})  is evaluated
by inserting approximations $\hat{u}_h$ and $\hat{v}_h,$ which are
tensor product B-splines as usual. 

To compute an entry in the stiffness matrix, 
we put (omitting the degrees)
$$\hat{u}_h(\xi,\eta)=N_i(\eta)M_j(\xi) , \quad \hat{v}_h(\xi,\eta)=N_k(\eta)M_{\ell}(\xi)$$
and insert these products into (\ref{num:weakwithspecF}).
For the first integral, it follows 
\begin{eqnarray}
\int_0^1 \int_0^1 
\frac{1}{J(\eta)}\xi \,\hat{u}_{h,\xi}(\xi,\eta)\hat{v}_{h,\xi}(\xi,\eta) \,\f{b}_1^T(\eta) \f{b}_1(\eta) \;\text{d}\xi\text{d}\eta \mbox{\hspace{5cm}}  \nonumber\\
=
\int_0^1{\xi \,M_j'(\xi)M_{\ell}'(\xi) \;\text{d}\xi}\cdot \int_0^1{\frac{1}{J(\eta)}N_i(\eta)N_k(\eta)\,\f{b}_1^T(\eta)\f{b}_1(\eta)\;\text{d}\eta}.\mbox{\hspace{1cm}}
\label{eq:twointegrals}
\end{eqnarray}
The two-dimensional integration can thus be carried out as the product of
two one-dimensional integrations, which is a great computational advantage
and which is an important consequence of the 
multiplicative structure of the Jacobian $\f{DF}$.
The other terms in the weak form (\ref{num:weakwithspecF}) possess
the same feature.  

The separation of variables that we observe in (\ref{eq:twointegrals}) is a special case of the so-called low rank tensor approximation that has
recently been introduced by Mantzaflaris et al.~\cite{Mantzaflaris2015,Mantzaflaris2017}.
In this approach, low rank approximations of the integral kernels 
are computed to provide a compact, separated representation of the integrals in IGA. In our case, there is no need to compute an approximation.
Instead, the parametrization provides directly a low rank tensor representation.

\subsection{Scaled Boundary IGA}\label{subs:scaledboundIGA}

In the SB-IGA, we transform first the PDE into the strong form 
(\ref{eq:strongxi}) in parametric coordinates, which leads to
\begin{equation}\label{eq:strongscaledIGA}
-\frac{\partial}{\partial \xi}\left(
{\small\frac{1}{J}}\left(\xi \f{b}_1^T\f{b}_1\hat{u}_{\xi}+\f{b}_2^T\f{b}_1\hat{u}_\eta\right)\right)-\frac{\partial}{\partial \eta}\left(
{\small\frac{1}{J}}\left( \f{b}_1^T\f{b}_2\hat{u}_{\xi}+\frac{1}{\xi} \f{b}_2^T\f{b}_2\hat{u}_\eta\right)\right)=\xi J \hat{f}.
\end{equation}
Then, a Galerkin projection with respect to the circumferential coordinate $\eta$ is derived, using the given representation of the boundary in terms of 
the B-splines $\f{N}(\eta)$. We  
set $\hat{u}(\xi,\eta)= \f{N}^T(\eta) \f{U}(\xi),\; \hat{v}(\xi,\eta)= \f{N}^T(\eta) \f{V},$ 
where $\f{U}(\xi) \in \mathbb{R}^n$ is the solution depending on the radial parameter $\xi$ and the variations $\f{V} \in \mathbb{R}^n$ are arbitrary.
We insert $\hat{u}$ in (\ref{eq:strongscaledIGA}), 
multiply by $\hat{v}$ and integrate with respect to $\eta$. This yields
\begin{align}\label{eq:multVintDeta}
-\int_0^1{\f{V}^T \f{N}{\small\frac{1}{J}}\f{b}_1^T \f{b}_1 \f{N}^T\f{U}_\xi \; \text{d}\eta } 
-\int_0^1{\f{V}^T \f{N}{\small\frac{1}{J}}\xi \,\f{b}_1^T \f{b}_1 \f{N}^T \f{U}_{\xi\xi} \; \text{d}\eta } \mbox{\hspace{3.2cm}} & 
\nonumber\\
-\int_0^1{\f{V}^T \f{N}{\small\frac{1}{J}}\f{b}_2^T \f{b}_1 \f{N}'^T\f{U}_\xi \; \text{d}\eta }
-\int_0^1{\f{V}^T \f{N}\frac{\partial}{\partial \eta }\left( {\small\frac{1}{J}}\left( \f{b}_1^T\f{b}_2\hat{u}_{\xi}+\frac{1}{\xi} \f{b}_2^T\f{b}_2\hat{u}_\eta\right)\right) \; \text{d}\eta  } & \nonumber \\
= \int_0^1{\f{V}^T \f{N}\,\xi J\hat{f} \; \text{d}\eta  }. &
\end{align}
The fourth integral of equation (\ref{eq:multVintDeta}) is integrated by parts, using the property $V_1 = V_n$ that holds for the variations of a periodic curve.
Thus it holds
\begin{align*}
-\int_0^1{\f{V}^T \f{N}\frac{\partial}{\partial \eta }\left( \frac{1}{J}\left( \f{b}_1^T\f{b}_2\hat{u}_{\xi} +\frac{1}{\xi} \f{b}_2^T\f{b}_2\hat{u}_\eta\right)\right) \; \text{d}\eta  } \mbox{\hspace{4cm}} \\
 =\int_0^1{\f{V}^T \f{N}'\frac{1}{J}\f{b}_1^T \f{b}_2\f{N}^T \f{U}_\xi \;\text{d}\eta}+
\int_0^1{\f{V}^T \f{N}'\frac{1}{J}\frac{1}{\xi}\f{b}_2^T \f{b}_2\f{N}^T \f{U} \;\text{d}\eta}.
\end{align*}
Since (\ref{eq:multVintDeta}) must hold for all variations $\f{V}$,
we can now perform the integration with respect to $\eta$ and generate a 
strong form for the vector $\f{U}(\xi)$. For this purpose,
we introduce a notation that is common in the SB-FEM and SB-IGA. 
Define $$\f{B}_1^T(\eta):=\frac{1}{J(\eta)}\f{N}\f{b}_1^T\in \mathbb{R}^{n\times 2},\quad \f{B}_2^T(\eta):=\frac{1}{J(\eta)}\f{N}'\f{b}_2^T\in \mathbb{R}^{n\times 2},$$ then (\ref{eq:multVintDeta}) becomes 
\begin{equation}\label{eq:ode-bvp}
\xi\f{M} \f{U}_{\xi\xi}+(\f{M}-\f{C}+\f{C}^T)\f{U}_\xi-\frac{1}{\xi}\f{K}\f{U}=\xi\f{S}(\xi),
\end{equation}
where
\begin{align*}
&\f{M}:=\int_0^1{\f{B}_1^T(\eta)\f{B}_1(\eta)J(\eta)\,\text{d}\eta},\\
&\f{C}:=\int_0^1{\f{B}_1^T(\eta)\f{B}_2(\eta)J(\eta)\,\text{d}\eta},\\
&\f{K}:=\int_0^1{\f{B}_2^T(\eta)\f{B}_2(\eta)J(\eta)\,\text{d}\eta},\\
&\f{S}(\xi):=-\int_0^1{\f{N}(\eta)J(\eta)\hat{f}(\xi,\eta)\,\text{d}\eta}.
\end{align*}
We obtain in this way an ODE in the scaling direction $\xi$, with boundary condition $\f{U}(\xi=1) = \f{0}$ and periodicity condition $U_1(\xi) = U_n(\xi)$.

There are several options for the numerical treatment of the ODE (\ref{eq:ode-bvp}) in scaling direction, as proposed in 
\cite{Chen2015,CHEN2016777,Natarajan2015}. Besides an analytical approach that is based on the solution of an eigenvalue problem with a 
Hamiltonian matrix, see Section 4 below, there is also the possibility of using collocation or standard Galerkin projection with respect to the $\xi$-variable.
The latter is based on a weak form in $\xi$ which possesses the same structure and the same separation property as shown in (\ref{eq:twointegrals}). 
For more details, we refer to the above references.

%% file: singularity.tex
\section{The Singularity in the Scaling Center}


In this section, we investigate the singularity in the scaling center
of the parametri\-zation in the case $d=2$ in more details.
It should be stressed that we do not face a singularity of the solution
of Poisson's equation (\ref{num:poisson}) here. On the contrary,
such an elliptic boundary value problem exhibits a regular behavior, with
the maximum principle bounding the solution of the strong form.
Coercivity and continuity of the bilinear form in the
weak formulation (\ref{num:weakform}), on the other hand, 
yield bounds and stability estimates in the energy and $H^1$-norms.

\subsection{The Boundary Value Problem in Scaling Direction}
To better understand the singularity at $\xi = 0$, we make use
of the boundary value problem that is generated by the SB-IGA.
In its original form, the SB-FEM, which provides the basic idea for SB-IGA,
is a semi-analytical method where the boundary value problem with respect to
$\xi$ is solved exactly \cite{Song2004}. We adopt here this approach since it offers insight
into the nature of the singularity. For a more general treatment of
singularities in IGA see \cite{TakacsJuettler2011,TakacsJuettler2012}.

We consider (\ref{eq:ode-bvp}) as homogeneous system
(writing $\f{U}_\xi = \f{U}^\prime$)
\begin{equation}\label{eq:odebvpnew}
\xi^2 \f{M} \f{U}^{\prime\prime}
                     + \xi (\f{M} - \f{C} + \f{C}^T) \f{U}^\prime 
                     - \f{K} \f{U} = \f{0}
\end{equation}
with a symmetric positive definite matrix $\f{M}  \in \mathbb{R}^{n\times n}$, a symmetric positive semi-definite matrix 
$\f{K}  \in \mathbb{R}^{n\times n}$, and 
a matrix $\f{C}  \in \mathbb{R}^{n\times n}$. 
Next, we introduce the new variables
\begin{equation}\label{eq:defy}
\f{y}(\xi) :=       \left( 
      \begin{array}{c} \f{U}(\xi) \\
       \f{W}(\xi)
      \end{array} \right)
      \quad \mbox{where} \,\,
       \f{W} 
      := \xi \f{M} \f{U}^\prime + \f{C}^T \f{U}.
\end{equation}
The second order ODE is thus transformed to a first order system
\begin{equation}\label{eq:hamiltonian}
\xi \f{y}^\prime = - \f{H} \f{y}
\end{equation}
with the Hamiltonian matrix
\[   \f{H} := \left( 
      \begin{array}{cc} \f{M}^{-1} \f{C}^T & - \f{M}^{-1}  \\
       - \f{K} + \f{C} \f{M}^{-1} \f{C}^T & - \f{C}\f{M}^{-1} 
      \end{array} \right)
\]
In other words, $\f{H} \in \mathbb{R}^{2n \times 2n}$ becomes a symmetric matrix when multiplied by the skew-symmetric matrix $\f{J}$,
\[ \f{J} = \left( 
      \begin{array}{cc} \f{0} &  \f{I}_n  \\
       - \f{I}_n & \f{0} 
      \end{array} \right) \quad \Rightarrow \,\, (\f{J} \f{A})^T =  \f{J} \f{A}.
\]
The characteristic polynomial 
$d(\lambda) = \mbox{det}(\lambda \f{I}_{2n} - \f{H})$ is an
even function, which means that the eigenvalues of $\f{H}$ 
come in pairs $(-\lambda_i, \lambda_i)$  
with $\mbox{Re} \lambda_i \geq 0$.
Hence the solution of (\ref{eq:hamiltonian}) in the homogeneous case is a linear combination of 
terms  
$c_i \xi^{\lambda_i} \fg{\phi}_i + \hat{c}_i \xi^{-\lambda_i}
\hat{\fg{\phi}_i}$ with corresponding eigenvectors $\fg{\phi}_i, 
\hat{\fg{\phi}}_i$.
Since the solution is finite in $\xi = 0$, 
one concludes that $\hat{c}_i = 0$, which 
cancels the instable part.

In practice, the numerical solution of the eigenvalue problem for $\f{H}$ 
is only feasible for comparably small dimensions. But our interest here lies 
on the insight that we obtain from it. The singularity in the scaling center 
thus looses its threat. However, the question remains what happens if we 
apply a discretization with respect to $\xi$ and do not utilize the 
eigenvalue solution.

\subsection{Practical Treatment of the Singularity}

We concentrate next on Galerkin-based IGA and discuss the 
practical treatment of the singularity in the scaling center.
In all numerical experiments that we performed so far, the singularity did not really influence the results, and we did not observe an instability nor a singular stiffness matrix. Can we explain this?

Consider the bilinear form (\ref{num:weakwithspecF}) in parametric coordinates. The fourth term in the integral is the critical one and it contains the factor $1/\xi$. Analogously to the separation in (\ref{eq:twointegrals}),
this term can also be written as the product of a well-defined integral with
respect to $\eta$ times the integral 
\begin{equation}\label{eq:integralsing}
\int_0^1\frac{1}{\xi} \,M_j(\xi)M_{\ell}(\xi) \;\text{d}\xi \,.
\end{equation}
We apply linear B-splines and analyze the contribution of the first element 
from $0$ to $h$ where $h$ is the mesh size. The integral (\ref{eq:integralsing}) then yields for $M_j=M_\ell=M_1$
\[    \int_0^h\frac{1}{\xi} (1-\xi/h) (1-\xi/h) \;\text{d}\xi 
 = \int_0^h   \left( \frac{1}{\xi}  - \frac{2}{h} + \frac{\xi}{h^2} \right) \;\text{d}\xi \,.
\]
Obviously, the integral over $1/\xi$ diverges while the other terms are 
not critical. In a numerical implementation, however, 
quadrature is used instead of exact integration. 
If we apply the midpoint rule as  simple example, we get 
\begin{equation}\label{eq:midpoint}
\int_0^h   \left( \frac{1}{\xi}  - \frac{2}{h} + \frac{\xi}{h^2} \right) \;\text{d}\xi \doteq 
h \left( \frac{1}{\xi}  - \frac{2}{h} + \frac{\xi}{h^2} \right)|_{\xi=h/2}
= \frac{1}{2} .
\end{equation}
The evaluation is thus independent of $h$,  and  the contribution of this integral to the stiffness matrix is always well-defined.

The above reasoning applies also to B-splines of higher degrees and to
higher order quadrature rules as long as the nodes of the quadrature rule are
in the interior of $[0, h ]$. In this way, we can conclude that the 
computation of the stiffness matrix in Galerkin-based IGA is not 
affected by the singularity, cf.~\cite{TakacsJuettler2011}.
However, at the moment we cannot say whether 
the convergence of the method is impaired or whether the condition number might become problematic.

\begin{figure}[h!]
	\centering
	\includegraphics[height = 3.95cm]{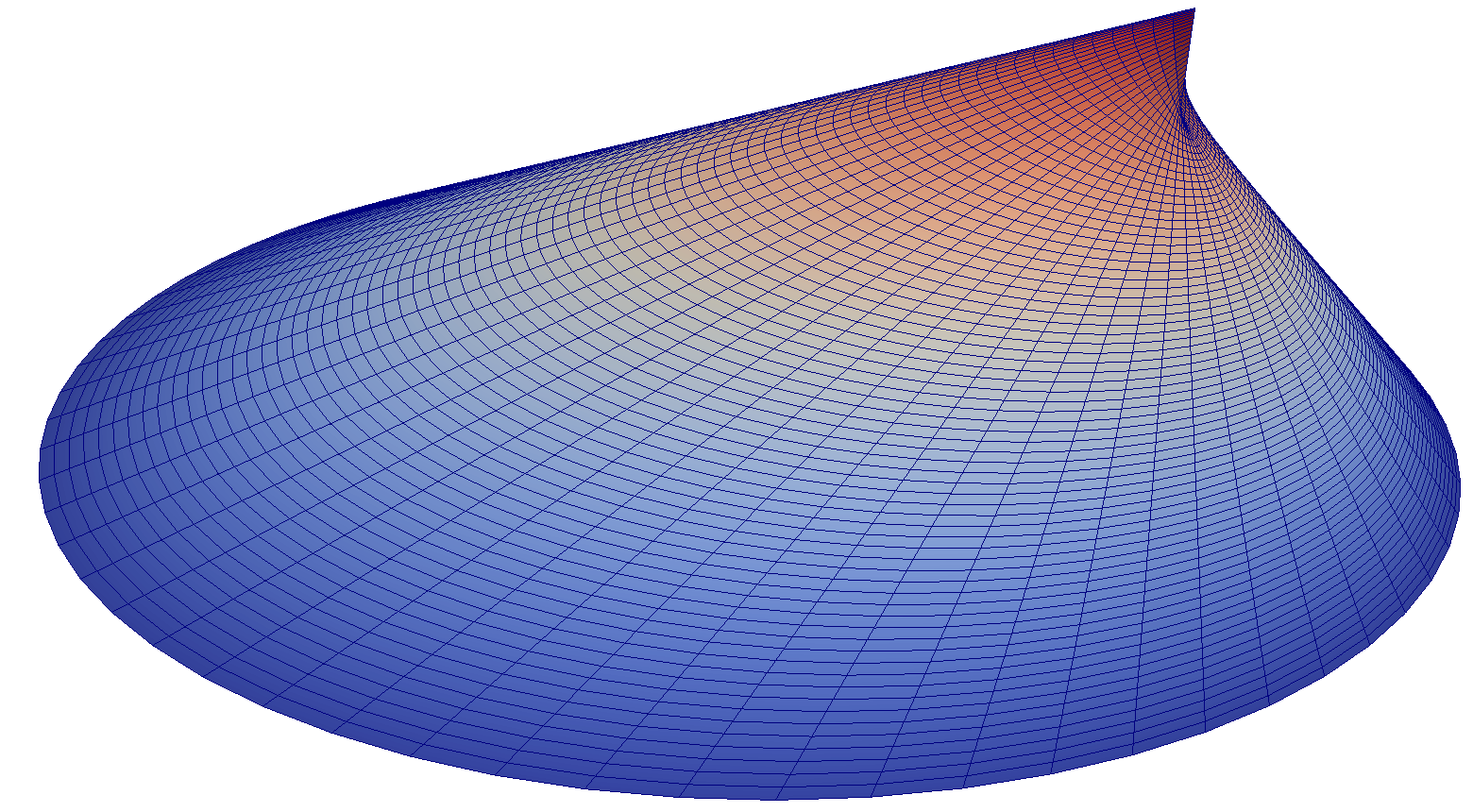}
	\includegraphics[height = 3.95cm]{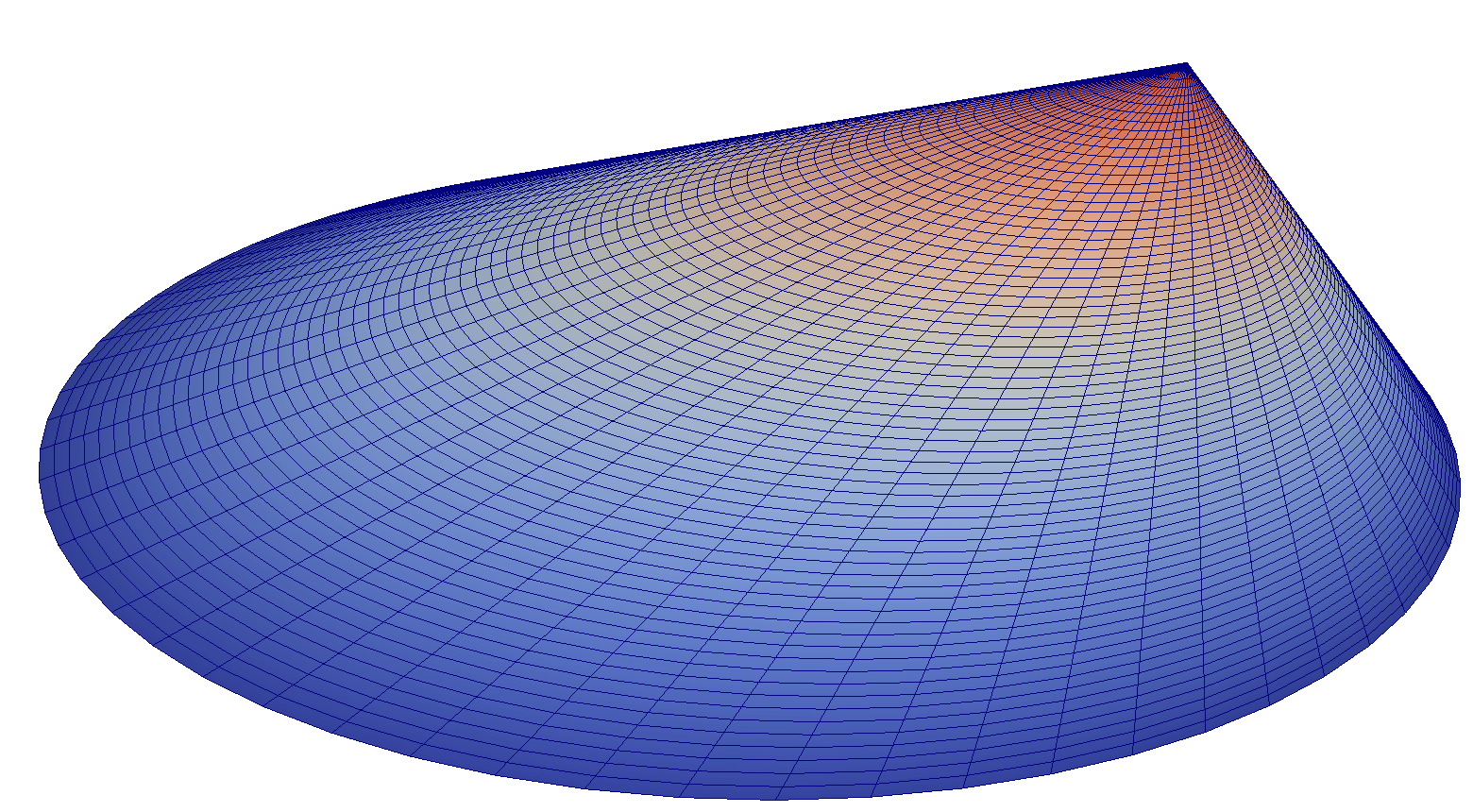}
	\caption{Off-center scaled boundary parametrization of a unit circle. Discontinuous (left) and constrained continuous (right) solutions. }
	\label{fig:singularity}
\end{figure}

Another problem is the smoothness of the numerical solution in the scaling center. The geometry mapping (\ref{def:geoFc}) projects the entire side of the parametric domain into the scaling center. Therefore several degrees of freedom correspond to the value there. By construction the solution is continuous around the singularity, however, only neighboring degrees of freedom are related to each other. As a consequence, their values across the singularity do not have to be the same, which results in the discontinuous numerical solution, see Fig.~\ref{fig:singularity}.

The continuity of the solution can be restored by imposing an equality constraint on all the degrees of freedom in the scaling center. However, this leads to a saddle point problem if Lagrange multipliers are used. It is also worth mentioning that the significant difference between the initially discontinuous and the constrained continuous solutions is only noticeable for very coarse meshes, while for the number of degrees of freedom large enough the solutions are nearly identical, see Fig.~\ref{fig:sqPlots} in Section 5.

%% file: computational.tex
\section{Numerical Examples}
We study two numerical examples to illustrate the usage of SB-parametrizations in IGA as described in Section 3, and we compare it with the standard rectangular parametrization. The simulations are performed using G+Smo, an open-source C++ library implementing main IGA routines \cite{jlmmz2014}.

\subsection{Poisson's Equation on a Unit Square}
\begin{figure}[h!]
	\centering
	\includegraphics[height = 3.55cm]{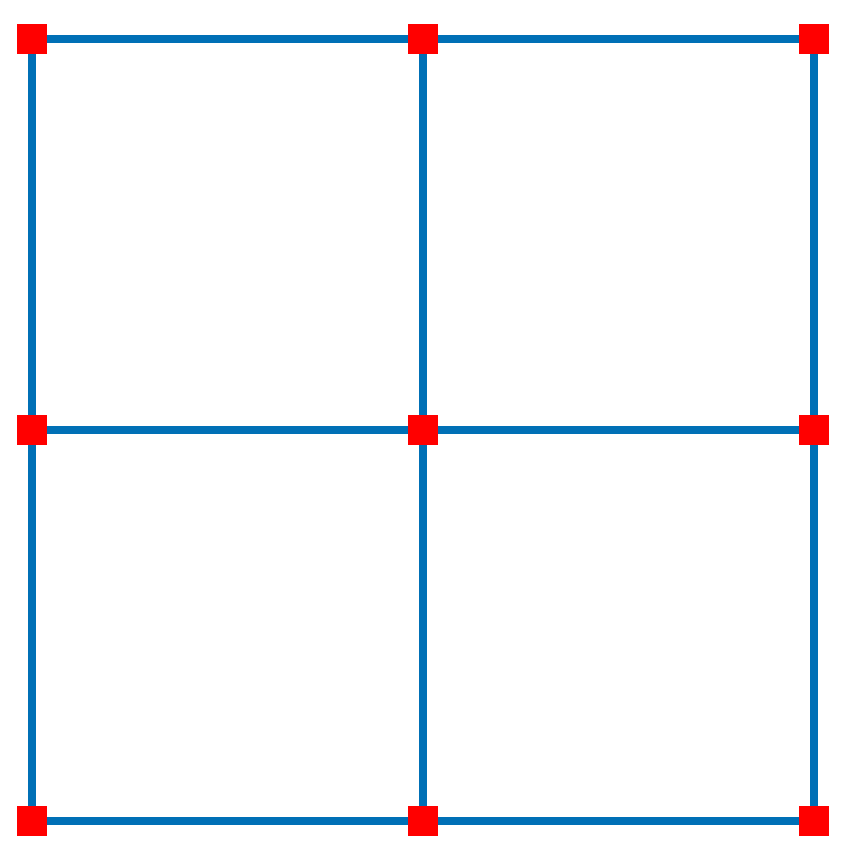}
	\includegraphics[height = 3.55cm]{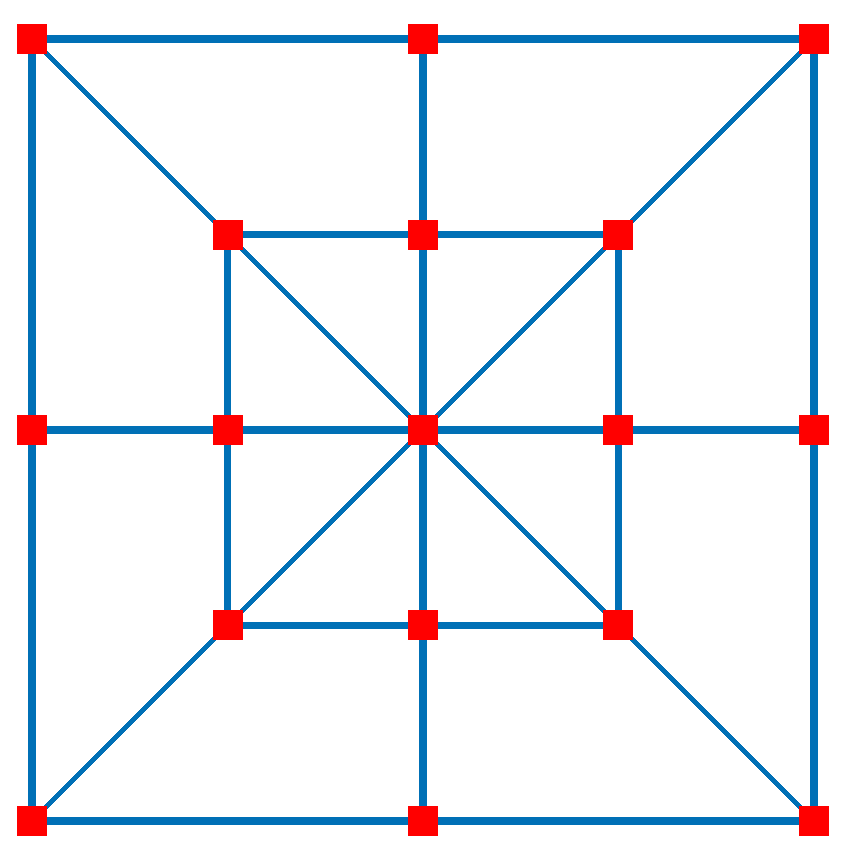}
	\includegraphics[height = 3.55cm]{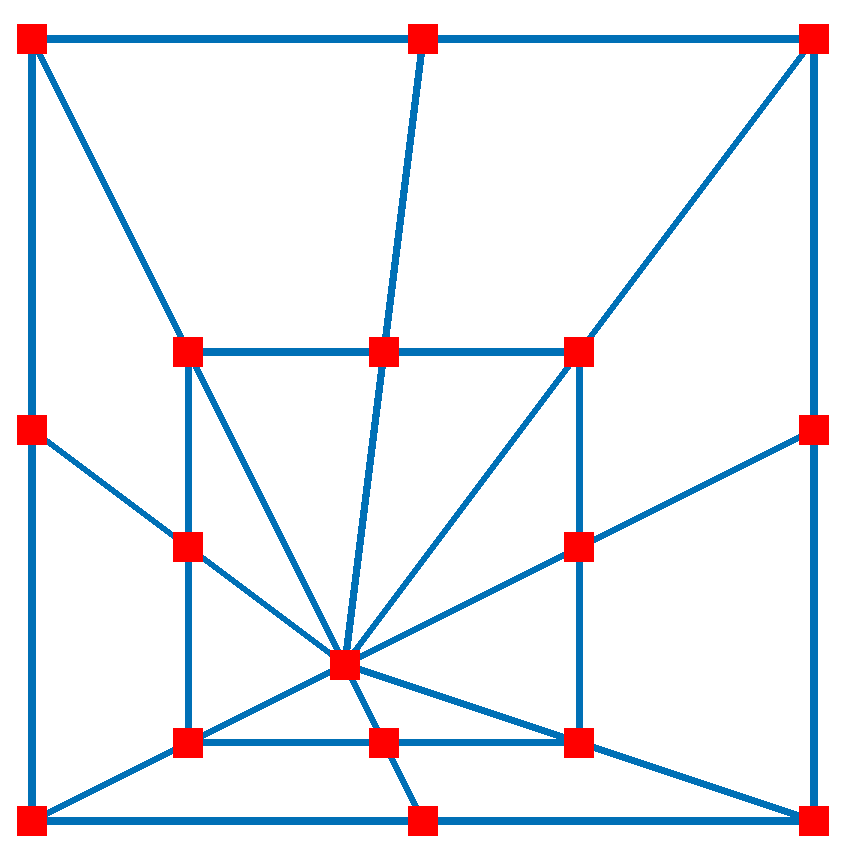}
	\includegraphics[height = 3.55cm]{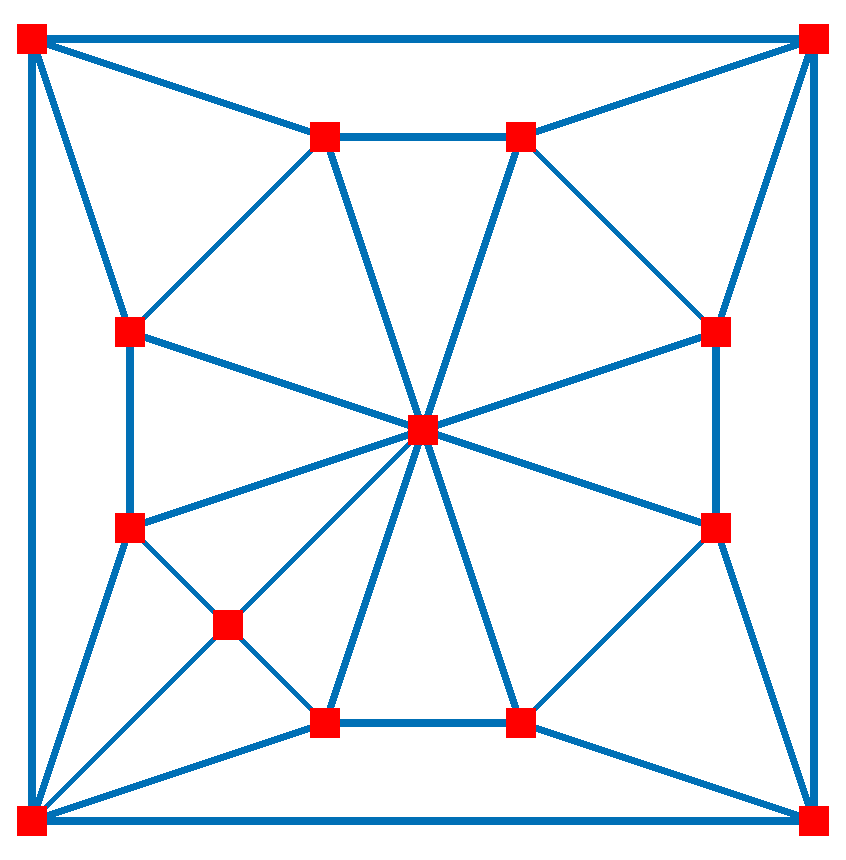}\\
	\includegraphics[height = 3.55cm]{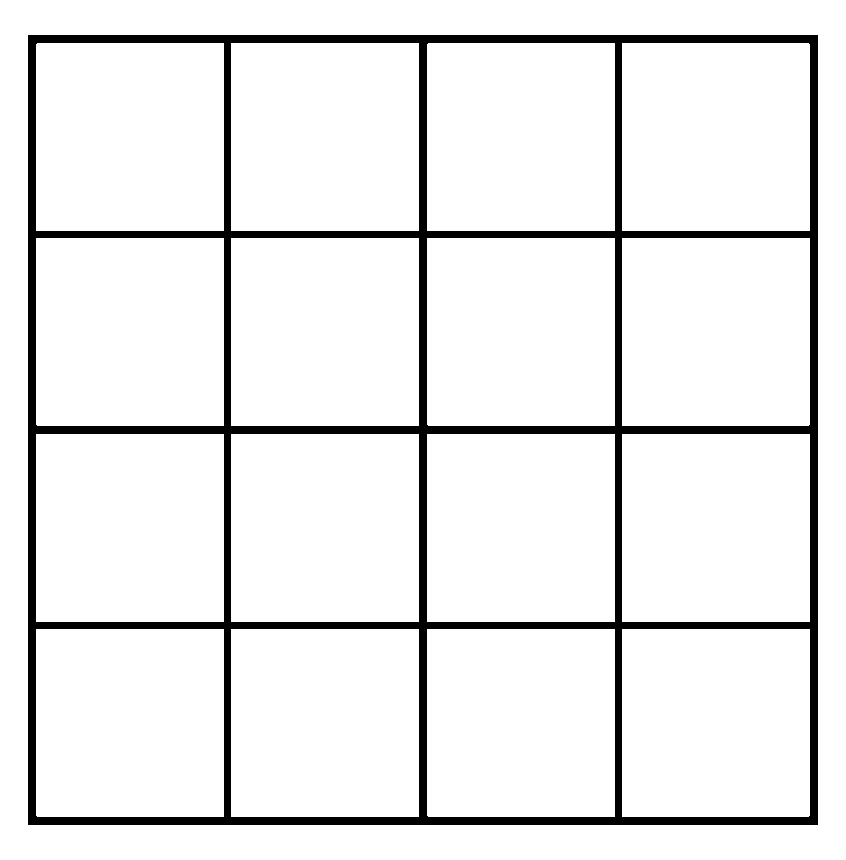}
	\includegraphics[height = 3.55cm]{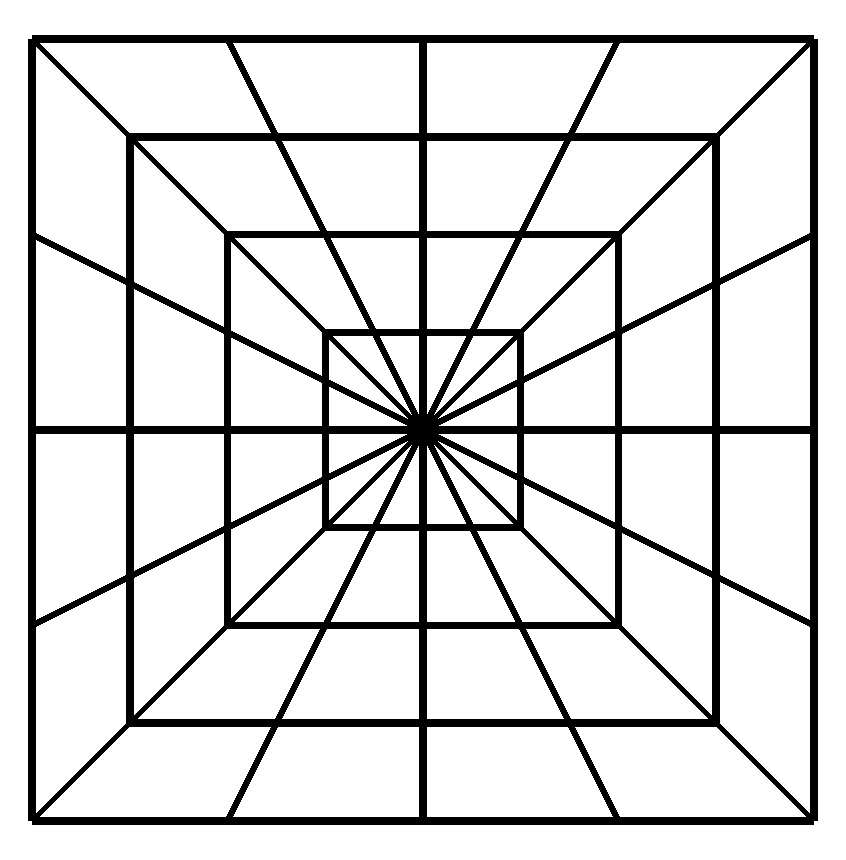}
	\includegraphics[height = 3.55cm]{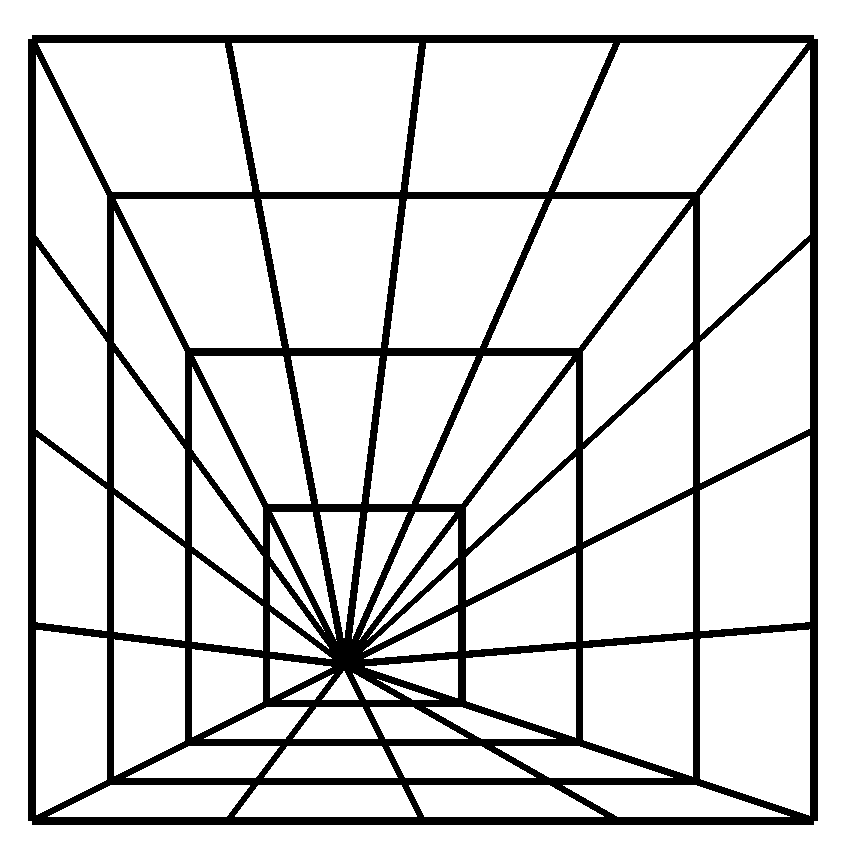}
	\includegraphics[height = 3.55cm]{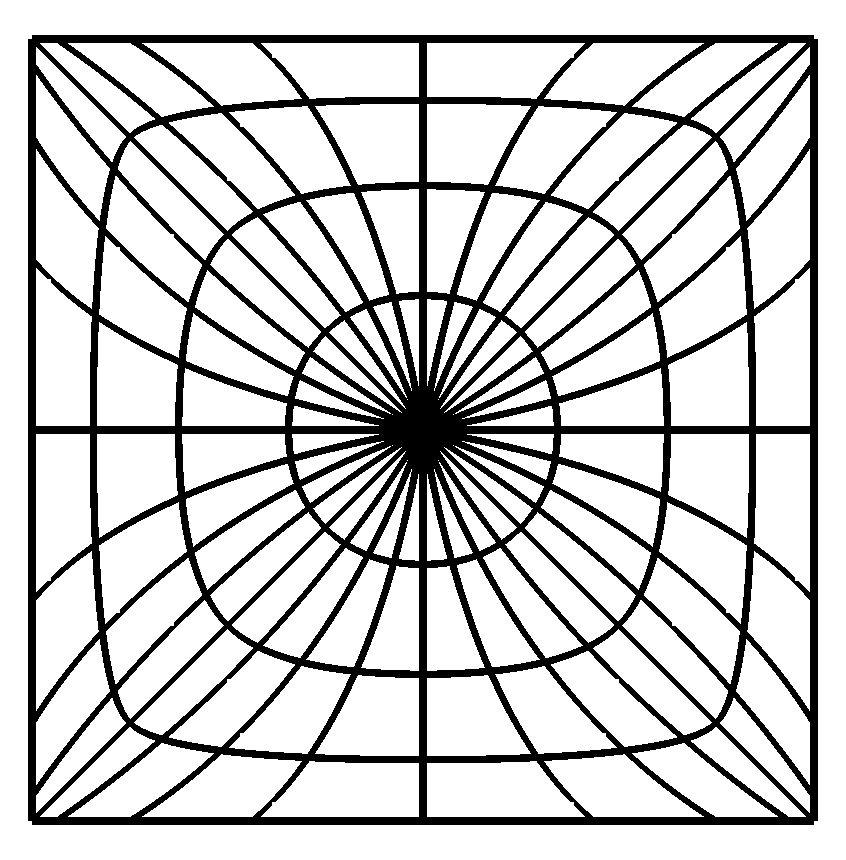}\\
	\includegraphics[height = 3.5cm]{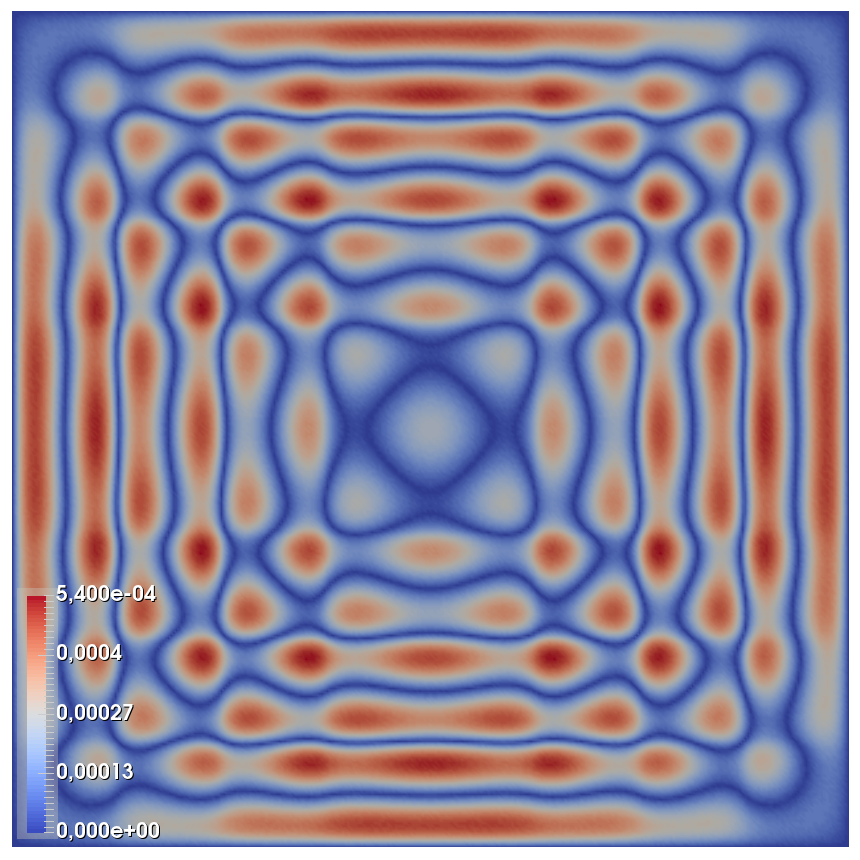}
	\includegraphics[height = 3.5cm]{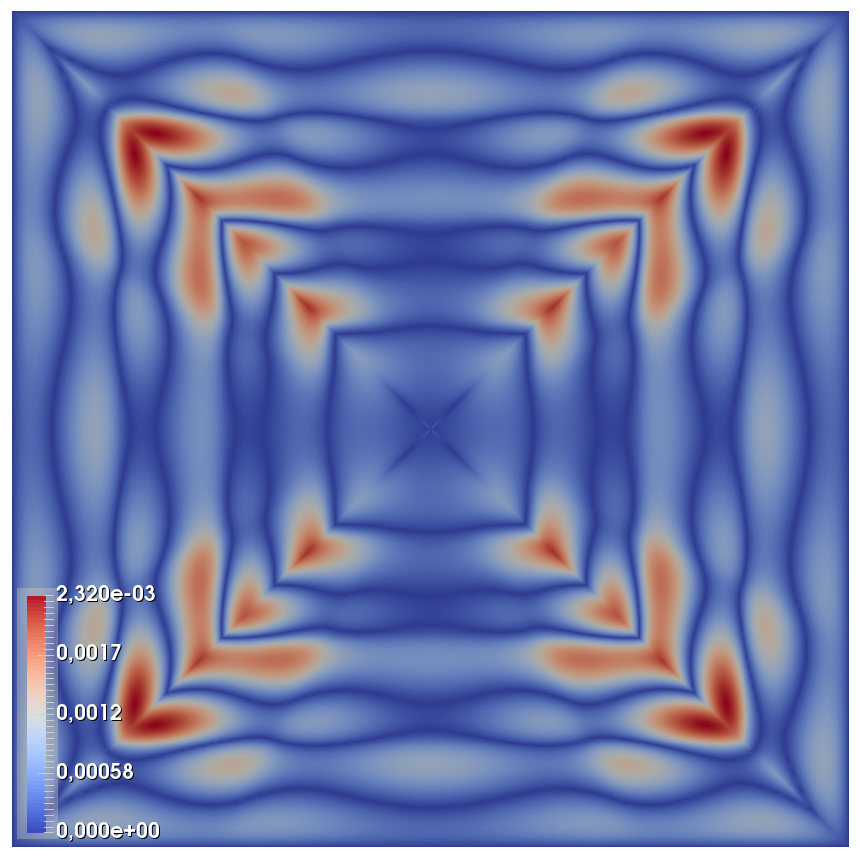}
	\includegraphics[height = 3.5cm]{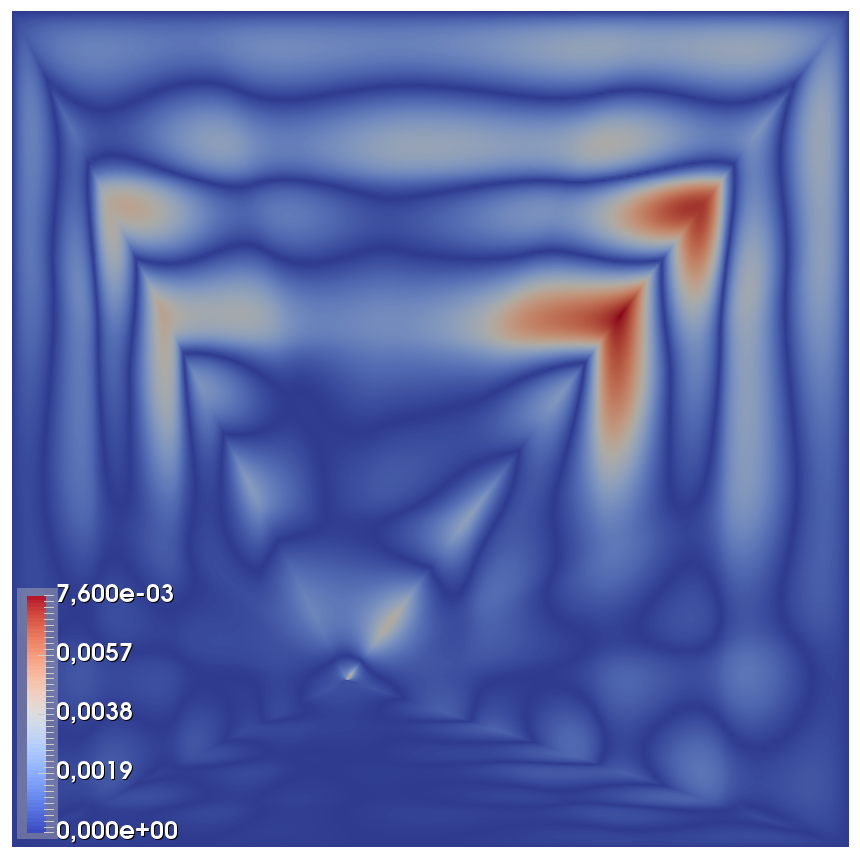}
	\includegraphics[height = 3.5cm]{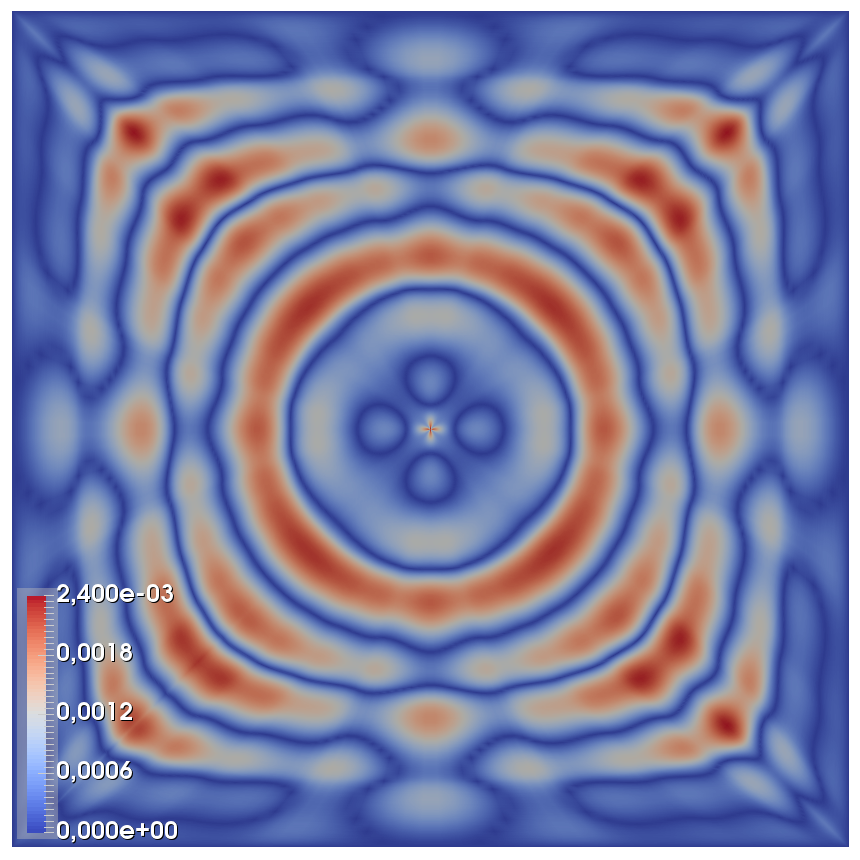}
	\caption{Unit square example. Rowwise: control nets with control points, isoparametric lines after two successive applications of h-refinement, numerical error for $\sim 10^2$ degrees of freedom. Columnwise: re\textbf{}ctangular, center scaled, off-center scaled and internally smooth parametrizations.}
	\label{fig:squares}
\end{figure}

We solve the Poisson problem on a unit square. The simplicity of this geometry allows a large variety of parametrizations to be generated. We select the following for comparison, see Fig.~\ref{fig:squares}: the {\it center scaled} as the most natural; the {\it off-center scaled} to investigate the influence of the scaling center location on the numerical error; the {\it internally smooth} to demonstrate that it is possible to generate an SB-parametrization that does not inherit non-smooth features from the boundary; and the {\it standard rectangular} to serve as a reference for comparison. Quadratic splines are used to represent the geometry and to discretize the equation for all parametrizations. We refer to Appendix A for the corresponding knot vectors and control points.

We consider the function 
\begin{equation}\label{eq:squareSol}
u^* = \cos\pi(x-0.5)\cos\pi(y-0.5),
\end{equation}
and we impose homogeneous boundary conditions $u=0$ on $\partial\Omega$. The function $u^*$ obviously fulfills the boundary conditions, thus it is the unique solution to (\ref{num:poisson}) with $f = -\Delta u^*$. We solve (\ref{num:poisson}) using different parametrizations and the known analytical solution is used to compute the numerical error, which is depicted in Fig.~\ref{fig:squares}.

After analyzing the error distribution it is immediately clear that the rectangular parametrization yields the best results. However, it should be noted for the SB-parametrizations that despite the singularity in the scaling center the numerical error is not observably higher there than in the other regions. The center and off-center scaled parametrizations indicate large error concentrations along the $C^0$ rays from the corners, whereas for the internally smooth parametrization the error is more distributed and is slightly lower. 

\begin{figure}[h!]
	\centering
	\includegraphics[height = 5.45cm]{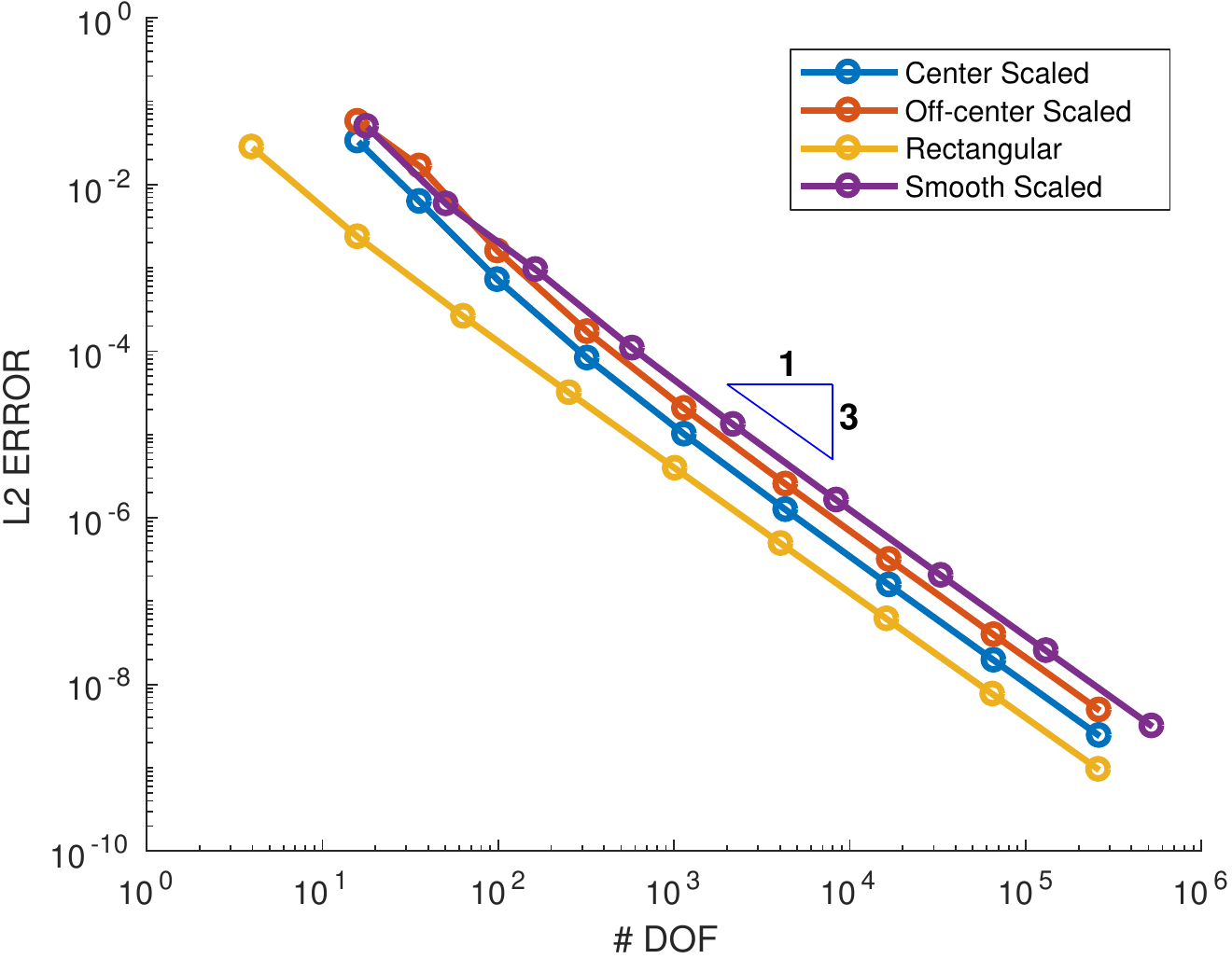}
	\includegraphics[height = 5.45cm]{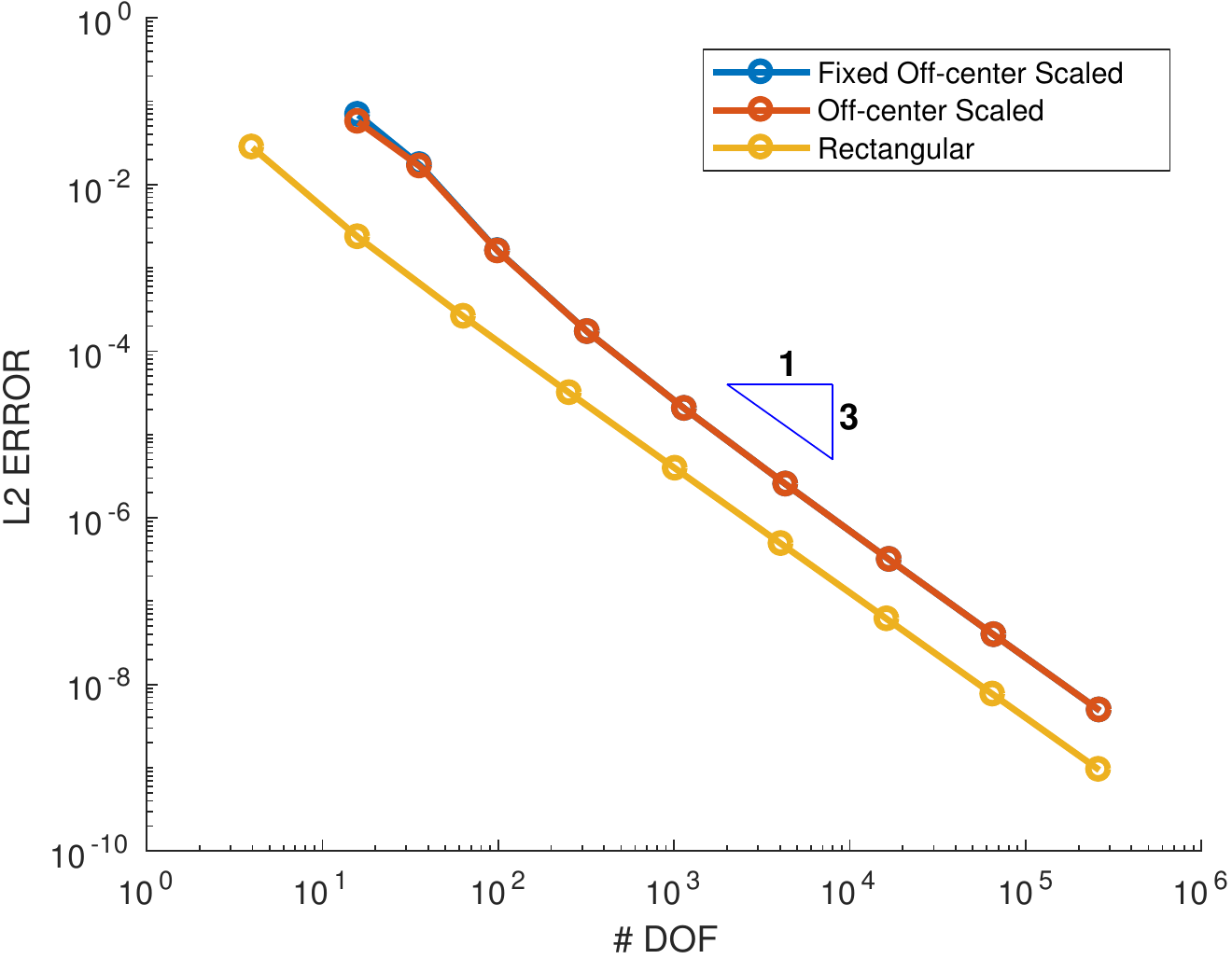}
	\caption{Unit square example, global error vs number of degrees of freedom. Left:  all parametrizations. Right: singularity treatment at the example of the off-center scaled parametrization.}
	\label{fig:sqPlots}
\end{figure}

We further compare the numerical error for different parametrizations by considering the global $L_2$ error while refining the analysis. The results are presented in Fig.~\ref{fig:sqPlots}. The convergence rate for all parametrizations is approximately of order three. The smallest error among the SB-parametrizations is found for the simplest and the most natural center scaled parametrization. Additionally, we do the convergence analysis for the off-center parametrization to study the influence of the singularity treatment as discussed in Subsection 4.2. The difference between the initially discontinuous and the constrained continuous solutions is only observed for a very coarse mesh, whereas for a sufficient number of degrees of freedom the two solutions are virtually identical. 

\subsection{Poisson's Equation on a Screw Compressor Rotor}
A screw compressor rotor with its intricate non-star-shaped geometry serves as an exemplary industrial application of SB-parametrizations due to a large variety of rotatory mechanisms in engineering. 

\begin{figure}[h!]
	\centering
	\includegraphics[height = 4cm]{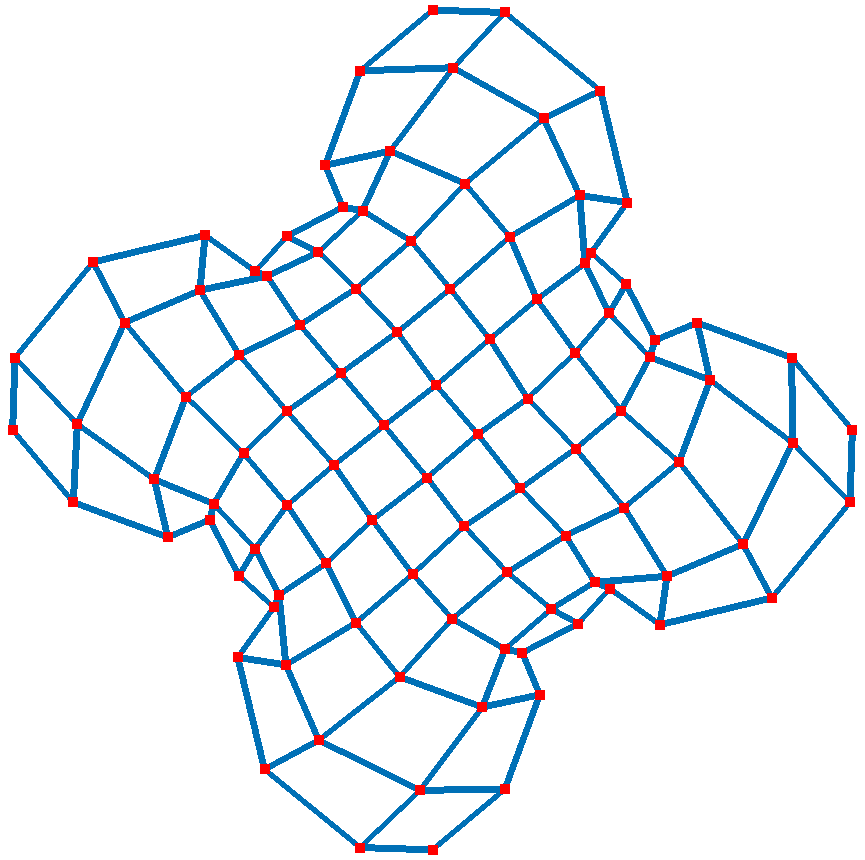}
	\includegraphics[height = 4cm]{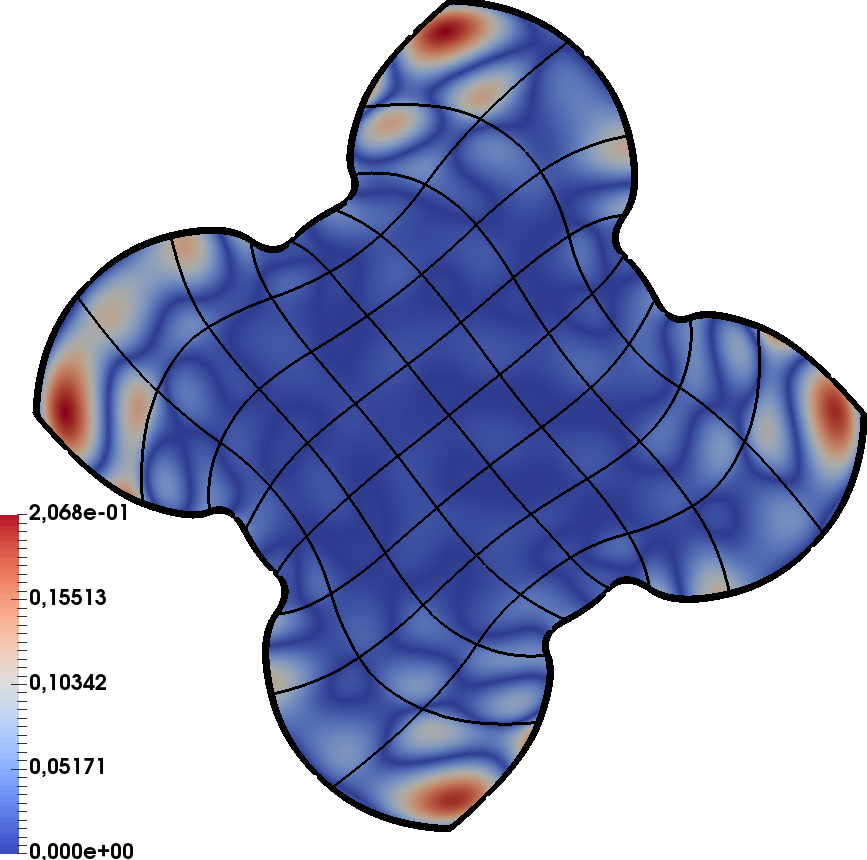}
	\includegraphics[height = 4cm]{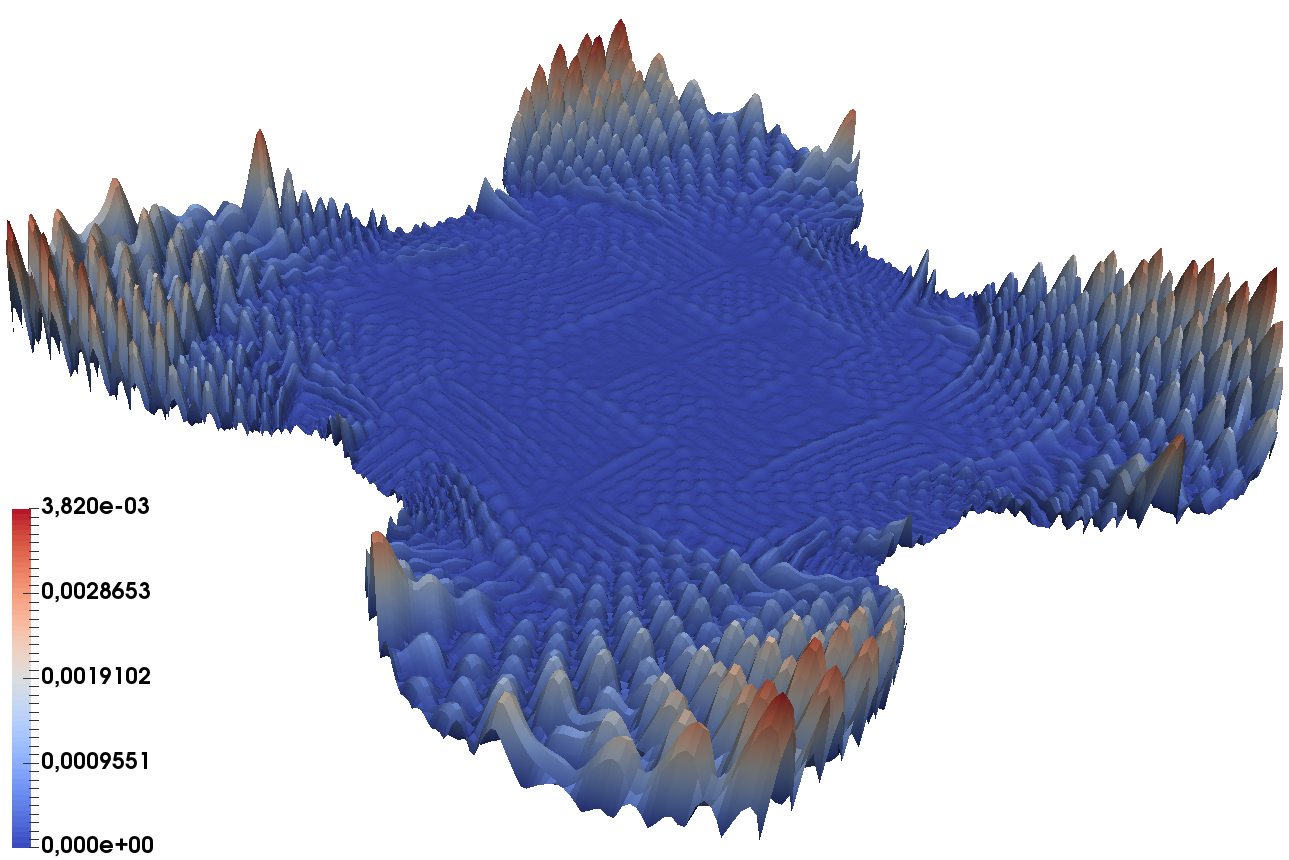}\\
	\includegraphics[height = 4cm]{SimFigures/cpRotor2.png}
	\includegraphics[height = 4cm]{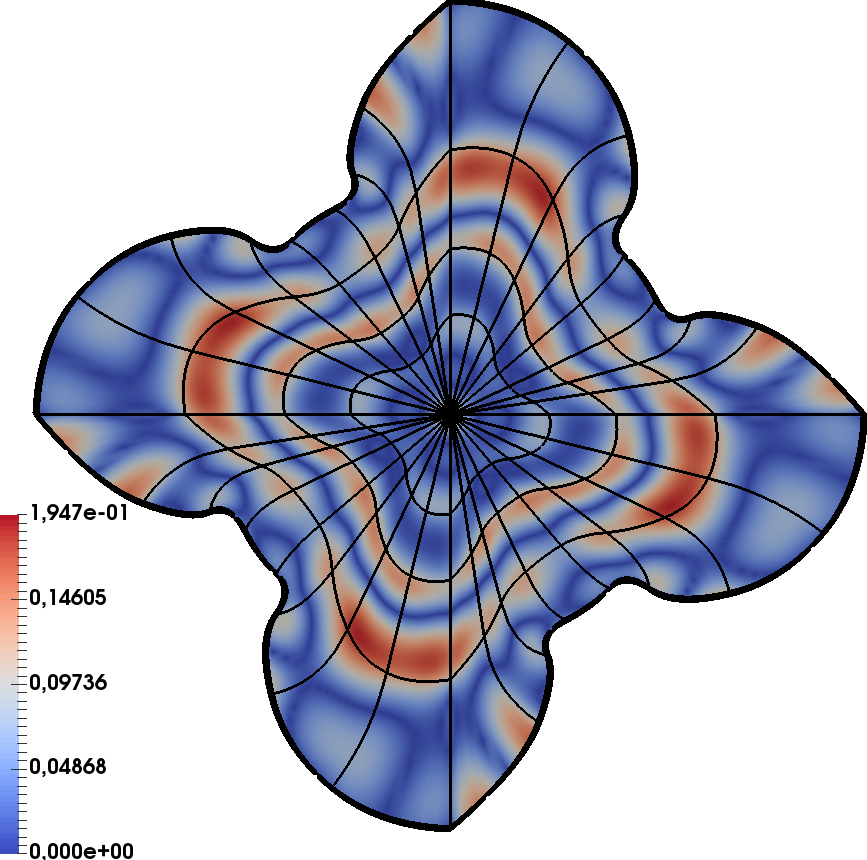}
	\includegraphics[height = 4cm]{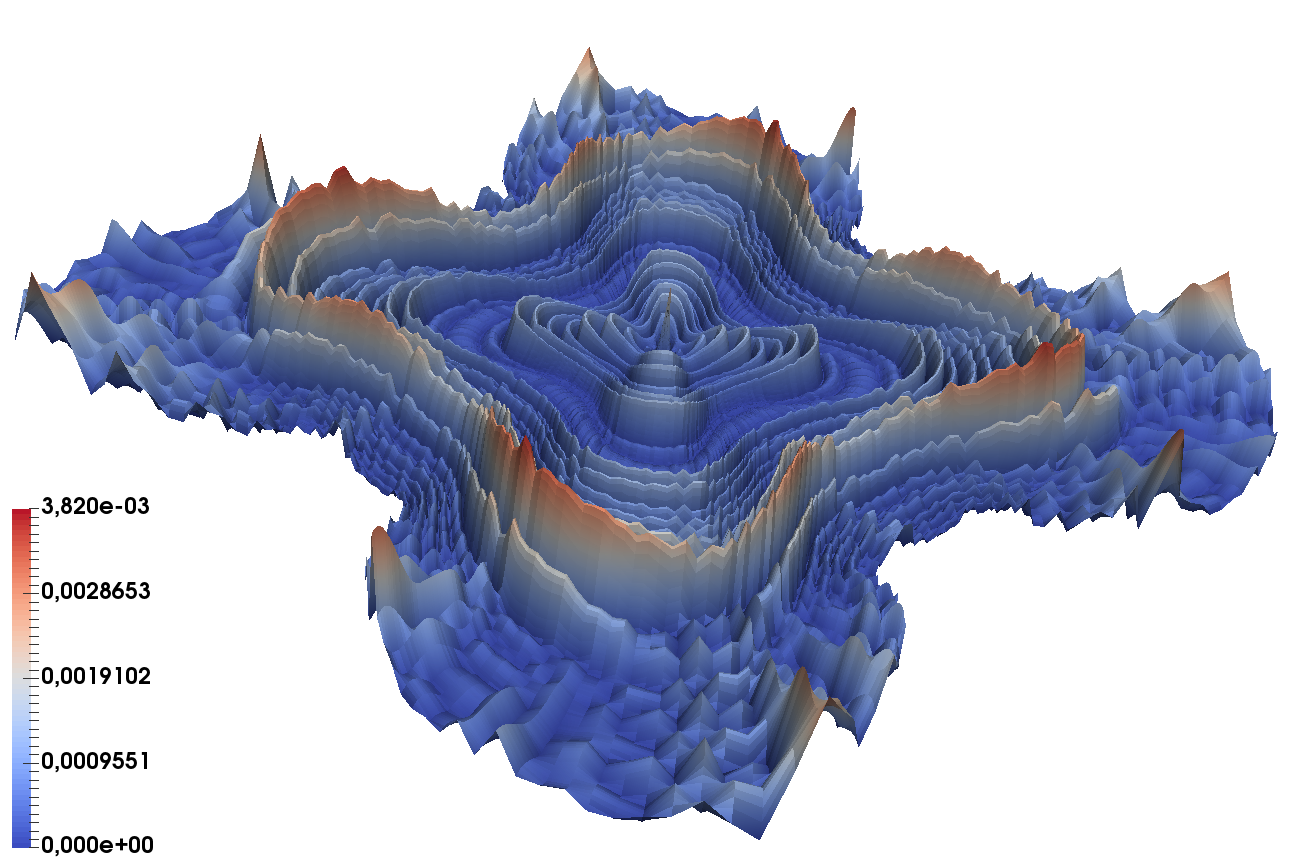}\\
	\caption{Rectangular and SB-parametrizations of a screw compressor rotor. Control nets with control points (left); isoparametric lines and numerical error for $\sim 10^2$ degrees of freedom (center); x1000 magnified numerical error as height field for $\sim 10^3$ degrees of freedom (right).}
	\label{fig:rotorsError}
\end{figure}

Here we consider the Poisson problem (\ref{num:poisson}) posed over the rotor geometry. We use the function
\begin{equation}\label{eq:rotorSol}
u^* = a^2 - x^2 - y^2 
\end{equation}
with the parameter $a \in \mathbb{R}$ as the Dirichlet boundary condition on $\partial\Omega$, which is enforced strongly. Then $u^*$ is the unique solution to (\ref{num:poisson}) with $f = -\Delta u^*$. The known analytical solution immediately allows to compute the numerical error, which can be used to measure the quality of the parametrization. The standard rectangular parametrization, which is relatively straightforward to generate in this case due to the presence of four sharp corners, serves as a reference for comparison, see Fig.~\ref{fig:rotorsError}.

\begin{figure}[h!]
	\centering
	\includegraphics[height = 6cm]{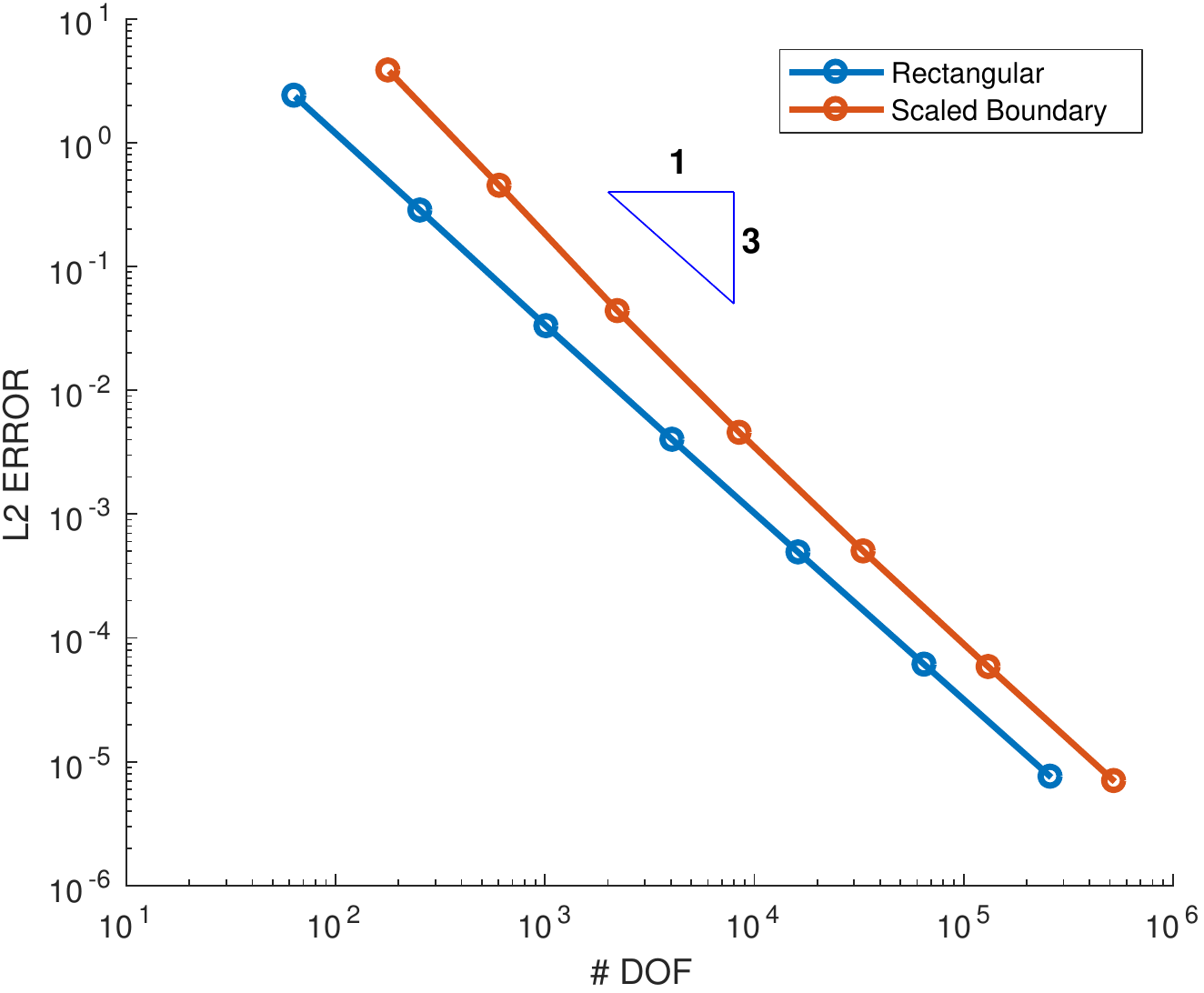}
	\caption{Screw compressor rotor. Global error vs number of degrees of freedom for rectangular and SB-parametrizations.}
	\label{fig:rotorsPlot}
\end{figure}

Probably the most remarkable result of this comparison is that the singularity in the scaling center does not give rise to any particularly relevant numerical error. Moreover, we note that both parametrizations yield errors of the same order of magnitude. To further study the numerical error, we use its $L_2$-norm as a global error and we compare the two parametrizations while refining the analysis. Fig.~\ref{fig:rotorsPlot} depicts the results. Both parametrizations demonstrate the convergence rate of roughly third order, and while the rectangular parametrizations offers slightly more accuracy per degree of freedom, the SB-parametrization exhibits a marginally higher convergence rate.

%% file: appendix.tex
\section*{Appendix A. Control Data for B-spline Objects}
\subsection*{A.1. Unit Square. Center Scaled Parameterization}
Here we describe the simplest and the most natural SB-parametrization of a unit square, which is shown in Fig.~\ref{fig:squares}. The basis is quadratic in both direction, with the knots vectors given by
\begin{equation*}
\Xi = \{0,0,0,0.25,0.25,0.5,0.5,0.75,0.75,1,1,1\}
\end{equation*}
and
\begin{equation*}
\Psi = \{0,0,0,1,1,1\}
\end{equation*}
and the control points from Table~\ref{table:center}.
The resulting planar parametrization inherits $C^0$ rays from the four corners of the square, remaining $C^1$ between them. 

\begin{table}[h]
	\caption{Unit square. Control points for the center scaled parametrization.}
	\label{table:center}
	\begin{tabularx}{\textwidth}{Y|YYYYYY}
		\hline
		$j$ & $\f{d}_{1,j}$   & $\f{d}_{2,j}$  & $\f{d}_{3,j}$   & $\f{d}_{4,j}$  & $\f{d}_{5,j}$   & $\f{d}_{6,j}$  \\
		& $\f{d}_{7,j}$   & $\f{d}_{8,j}$  & $\f{d}_{9,j}$   &            &             &            \\ \hline
		1   & (0,0)       & (0.5,0)    & (1,0)       & (1,0.5)    & (1,1)       & (0.5,1)    \\
		& (0,1)       & (0,0.5)    & (0,0)       &            &             &            \\
		2   & (0.25,0.25) & (0.5,0.25) & (0.75,0.25) & (0.75,0.5) & (0.75,0.75) & (0.5,0.75) \\
		& (0.25,0.75) & (0.25,0.5) & (0.25,0.25) &            &             &            \\
		3   & (0.5,0.5)   & (0.5,0.5)  & (0.5,0.5)   & (0.5,0.5)  & (0.5,0.5)   & (0.5,0.5)  \\
		& (0.5,0.5)   & (0.5,0.5)  & (0.5,0.5)   &            &             &            \\ \hline
	\end{tabularx}
\end{table}

\subsection*{A.2. Unit Square. Internally Smooth Scaled Parameterization}
The internally smooth parametrization of a unit square shown in Fig.~\ref{fig:squares} is created by adding additional control points at the corners. By doing so we preserve sharp $C^0$ corners, while eliminating three $C^0$ rays connecting them to the scaling center. The resulting parametrization is smooth everywhere except for the four corners and one unavoidable $C^0$ ray, where the periodic boundary conditions are to be imposed. The quadratic basis is used, the control points are given is Table~\ref{table:smooth} and the knots vector are
\begin{equation*}
	\Xi = \{0,0,0,0.125,0.25,0.375,0.5,0.625,0.75,0.875,1,1,1\}
\end{equation*}
and
\begin{equation*}
	\Psi = \{0,0,0,1,1,1\}.
\end{equation*}
Note that there is no knot repetition except for the end knots. Thus the internal parametrization is smooth. 
\begin{table}[h]
	\caption{Unit square. Control points for the internally smooth parametrization.}
	\label{table:smooth}
	\begin{tabularx}{\textwidth}{Y|YYYYY}
		\hline
		$j$ & $\f{d}_{1,j}$     & $\f{d}_{2,j}$     & $\f{d}_{3,j}$     & $\f{d}_{4,j}$     & $\f{d}_{5,j}$     \\
		& $\f{d}_{6,j}$     & $\f{d}_{7,j}$     & $\f{d}_{8,j}$     & $\f{d}_{9,j}$     & $\f{d}_{10,j}$    \\ \hline
		1   & (0,0)         & (0,0)         & (1,0)         & (1,0)         & (1,1)         \\
		& (1,1)         & (0,1)         & (0,1)         & (0,0)         & (0,0)         \\
		2   & (0.25,0.25)   & (0.375,0.125) & (0.625,0.125) & (0.875,0.375) & (0.875,0.625) \\
		& (0.625,0.875) & (0.375,0.875) & (0.125,0.625) & (0.125,0.375) & (0.25,0.25)   \\
		3   & (0.5,0.5)     & (0.5,0.5)     & (0.5,0.5)     & (0.5,0.5)     & (0.5,0.5)     \\
		& (0.5,0.5)     & (0.5,0.5)     & (0.5,0.5)     & (0.5,0.5)     & (0.5,0.5)     \\ \hline
	\end{tabularx}
\end{table}

%% file: ms.bbl
\begin{thebibliography}{10}

\bibitem{Berger2007}
M.~Berger.
\newblock {\em A Panoramic View of Riemannian Geometry}.
\newblock Springer, 2007.

\bibitem{Chen2015}
L.~Chen, W.~Dornisch, and S.~Klinkel.
\newblock Hybrid collocation-{G}alerkin approach for the analysis of surface
  represented 3{D}-solids employing {SB-FEM}.
\newblock {\em Comput. Meth. Appl. Mech. Engrg.}, 295:268--289, 2015.

\bibitem{CHEN2016777}
L.~Chen, B.~Simeon, and S.~Klinkel.
\newblock A {NURBS} based {G}alerkin approach for the analysis of solids in
  boundary representation.
\newblock {\em Computer Methods in Applied Mechanics and Engineering},
  305:777--805, 2016.

\bibitem{Cottrell.2009}
J.~A. Cottrell, T.~JR Hughes, and Y.~Bazilevs.
\newblock {\em Isogeometric analysis: toward integration of {CAD} and {FEA}}.
\newblock John Wiley \& Sons, 2009.

\bibitem{Boor1978}
C.~de~Boor.
\newblock {\em A Practical Guide to Splines}.
\newblock Springer, 1978.

\bibitem{DEDE2015807}
L.~Dedè and A.~Quarteroni.
\newblock Isogeometric analysis for second order partial differential equations
  on surfaces.
\newblock {\em Computer Methods in Applied Mechanics and Engineering}, 284:807
  -- 834, 2015.
\newblock Isogeometric Analysis Special Issue.

\bibitem{giannelli2012thb}
C.~Giannelli, B.~J\"uttler, and H.~Speleers.
\newblock {THB}-splines: the truncated basis for hierarchical splines.
\newblock {\em Computer Aided Geometric Design}, 29(7):485--498, 2012.

\bibitem{gravesen2014}
J.~Gravesen, A.~Evgrafov, D.-M. Nguyen, and P.~N{\o}rtoft.
\newblock Planar parametrization in isogeometric analysis.
\newblock In M.~Floater, T.~Lyche, M.-L. Mazure, K.~M{\o}rken, and L.~L.
  Schumaker, editors, {\em Mathematical Methods for Curves and Surfaces: 8th
  International Conference, MMCS 2012, Oslo, Norway, June 28 -- July 3, 2012,
  Revised Selected Papers}, pages 189--212. Springer Berlin Heidelberg, 2014.

\bibitem{Hughes2005}
T.~J.~R. Hughes, J.~A. Cottrell, and Y.~Bazilevs.
\newblock {I}sogeometric analysis: {C}{A}{D}, finite elements, {N}{U}{R}{B}{S},
  exact geometry and mesh refinement.
\newblock {\em Comput. Meth. Appl. Mech. Engrg.}, 194:4135--4195, 2005.

\bibitem{jlmmz2014}
B.~J\"uttler, U.~Langer, A.~Mantzaflaris, S.~Moore, and W.~Zulehner.
\newblock Geometry + simulation modules: Implementing isogeometric analysis.
\newblock {\em Proc. Appl. Math. Mech.}, 14(1):961--962, 2014.
\newblock Special Issue: 85th Annual Meeting of the Int. Assoc. of Appl. Math.
  and Mech. (GAMM), Erlangen 2014.

\bibitem{Klinkel2015}
S.~Klinkel, L.~Chen, and W.~Dornisch.
\newblock A {NURBS} based hybrid collocation-{G}alerkin method for the analysis
  of boundary represented solids.
\newblock {\em Comput. Meth. Appl. Mech. Engrg.}, 284:689--711, 2015.

\bibitem{Mantzaflaris2015}
A.~Mantzaflaris, B.~J\"uttler, B.N. Khoromskij, and U.~Langer.
\newblock {M}atrix generation in isogeometric analysis by low rank tensor
  approximation.
\newblock In J.-D. Boisonnat, A.~Cohen, O.~Gibaru, C.~Gout, T.~Lyche, M.-L.
  Mazure, and L.~L. Schumaker, editors, {\em Curves and Surfaces,}, volume
  9213, pages 321--340. LNCS, Springer, 2015.

\bibitem{Mantzaflaris2017}
A.~Mantzaflaris, B.~Jüttler, B.~N. Khoromskij, and U.~Langer.
\newblock Low rank tensor methods in {G}alerkin-based isogeometric analysis.
\newblock {\em Computer Methods in Applied Mechanics and Engineering},
  316:1062--1085, 2017.

\bibitem{Natarajan2015}
S.~Natarajan, J.~C. Wang, C.~Song, and C.~Birk.
\newblock Isogeometric analysis enhanced by the scaled boundary finite element
  method.
\newblock {\em Comput. Meth. Appl. Mech. Engrg.}, 283:733--762, 2015.

\bibitem{ORourke1987}
J.~O'~Rourke.
\newblock {\em Art Gallery Theorems and Algorithms}.
\newblock Oxford University Press, 1987.

\bibitem{Piegl1997}
L.~Piegl and W.~Tiller.
\newblock {\em The {NURBS} book}.
\newblock Monographs in visual communications. Springer, 1997.

\bibitem{SAPUTRA2015213}
A.~A. Saputra, C.~Birk, and C.~Song.
\newblock Computation of three-dimensional fracture parameters at interface
  cracks and notches by the scaled boundary finite element method.
\newblock {\em Engineering Fracture Mechanics}, 148:213 -- 242, 2015.

\bibitem{Song2004}
C.~Song.
\newblock A matrix function solution for the scaled boundary finite-element
  equation in statics.
\newblock {\em Comput. Methods Appl. Mech. Engrg.}, 193(23):2325--2356, 2004.

\bibitem{Song1997}
C.~Song and J.~P. Wolf.
\newblock The scaled boundary finite-element method-alias consistent
  infinitesimal finite-element cell method-for elastodynamics.
\newblock {\em Comput. Methods Appl. Mech. Engrg.}, 147:329--355, 1997.

\bibitem{Song2000}
C.~Song and J.~P. Wolf.
\newblock The scaled boundary finite-element method--a primer: solution
  procedures.
\newblock {\em Comput. Struct.}, 78(1):211--225, 2000.

\bibitem{TakacsJuettler2011}
T.~Takacs and B.~J{ü}ttler.
\newblock Existence of stiffness matrix integrals for singularly parameterized
  domains in isogeometric analysis.
\newblock {\em Computer Methods in Applied Mechanics and Engineering},
  200(49-52):3568--3582, 2011.

\bibitem{TakacsJuettler2012}
T.~Takacs and B.~J{ü}ttler.
\newblock H2 regularity properties of singular parameterizations in
  isogeometric analysis.
\newblock {\em Graphical Models}, 74(6):361--372, 2012.

\bibitem{Toshniwal20171005}
D.~Toshniwal, H.~Speleers, R.~R. Hiemstra, and T.~JR Hughes.
\newblock Multi-degree smooth polar splines: A framework for geometric modeling
  and isogeometric analysis.
\newblock {\em Computer Methods in Applied Mechanics and Engineering}, 316:1005
  -- 1061, 2017.
\newblock Special Issue on Isogeometric Analysis: Progress and Challenges.

\bibitem{utri2017energy}
M.~Utri and A.~Br{\"u}mmer.
\newblock Energy potential of dual lead rotors for twin screw compressors.
\newblock In {\em IOP Conference Series: Materials Science and Engineering},
  volume 232, page 012018. IOP Publishing, 2017.

\bibitem{vuong2011hierarchical}
A.-V. Vuong, C.~Giannelli, B.~J{\"u}ttler, and B.~Simeon.
\newblock A hierarchical approach to adaptive local refinement in isogeometric
  analysis.
\newblock {\em Computer Methods in Applied Mechanics and Engineering},
  200(49):3554--3567, 2011.

\bibitem{Xu2011}
G.~Xu, B.~Mourrain, R.~Duvigneau, and A.~Galligo.
\newblock Parametrization of computational domain in isogeometric analysis:
  methods and comparison.
\newblock {\em Computer Methods in Applied Mechanics and Engineering},
  200(23--24):2021--2031, 2011.

\end{thebibliography}
